\documentclass[12pt]{article}%
\usepackage{amsmath}
\usepackage{amsmath}
\usepackage{amssymb}
\usepackage{dsfont}
\usepackage{cite}

\usepackage{amsthm}

\newcommand{\eT}{{\mathbf{T}}}
\def\nn{\nonumber}

\def\A{\mathfrak A}

\def\b{\beta}

\def\g{\gamma}
\def\G{\Gamma}
\def\de{\delta}

\def\Lam{\Lambda}

\def\la{\lambda}
\def\om{\omega}

\def\th{\theta}
\def\Th{\Theta}
\def\Ups{\Upsilon}

\def\vp{\varphi}
\def\vt{\vartheta}
\def\ve{\varepsilon}
\def\wh{\widehat}
\def\wt{\widetilde}
\def\ov{\overline}
\def\br{\breve}

\def\BC{{\mathbb C}}

\def\BR{{\mathbb R}}
\def\BN{{\mathbb N}}
\def\BZ{{\mathbb Z}}
\def\clp{{\mathcal P}}
\def\cla{{\mathcal A}}

\def\cle{{\mathcal E}}

\def\clh{{\mathcal H}}

\def\cli{{\mathcal I}}
\def\clk{{\mathcal K}}

\def\cln{{\mathcal N}}
\def\clu{{\mathcal U}}

\def\cld{{\mathcal D}}

\def\spa{{\rm Span}}
\def\diag{\mathrm{diag}}
\newcommand{\E}{\mathrm{e}}
\newcommand{\I}{\mathrm{i}}

\newtheorem{Pa}{Paper}[section]
\newtheorem{Tm}[Pa]{{\bf Theorem}}
\newtheorem{La}[Pa]{{\bf Lemma}}
\newtheorem{Cy}[Pa]{{\bf Corollary}}
\newtheorem{Rk}[Pa]{{\bf Remark}}

\newtheorem{Dn}[Pa]{{\bf Definition}}
\newtheorem{Nn}[Pa]{{\bf Notation}}
\newtheorem{Pn}[Pa]{{\bf Proposition}}

\title{Dirac systems with locally 
square-integrable potentials: direct and inverse problems \\ for the spectral functions}

\author{Alexander Sakhnovich}

\date{}
\begin{document}
\maketitle

\begin{flushright}
Faculty of Mathematics,
University
of
Vienna, \\
Oskar-Morgenstern-Platz 1, A-1090 Vienna,
Austria, \\
e-mail: {\tt oleksandr.sakhnovych@univie.ac.at}
\end{flushright}

\vspace{0.3em}

\begin{flushleft}
{\it Dedicated to Fritz Gesztesy on the occasion of his 70th birthday}
\end{flushleft}

\begin{abstract}  
We solve the inverse problems to recover Dirac systems on an interval or semiaxis from their 
spectral functions (matrix valued functions) for the case of locally square-integrable potentials.
Direct problems in terms of spectral functions are treated as well.
Moreover, we present necessary and sufficient conditions on the given distribution matrix valued function to be
a spectral function of some Dirac system with a locally square-integrable potential.
Interesting connections with Paley-Wiener sampling measures appear in the case of scalar
spectral functions.
\end{abstract}

\vspace{0.3em}

{MSC(2020): 34L40, 34A55, 47B15}

\vspace{0.2em}

{\bf Keywords}: Dirac system, locally square-integrable potential, spectral matrix function, inverse problem, characterisation of spectral functions, Paley-Wiener sampling measure. 

\section{Introduction} \label{intro}
\setcounter{equation}{0}
{\bf 1.} Dirac (Dirac-type) system is a classical object of analysis. It is also called ZS or AKNS, which signifies
its various applications. This system may be written down in the form
\begin{align} &       \label{1.1}
\frac{d}{dx}y(x, z )=\I (z j+jV(x))y(x,
z ) \quad
(x \geq 0),
\end{align} 
where $\I$ is the imaginary unit $(\I^2=-1)$, $z$ is a complex-valued spectral parameter, $j=\diag\{I_{m_1},-I_{m_2}\}$ (i.e., $j$ is a block diagonal matrix with the
blocks $I_{m_1}$ and $-I_{m_2}$ on the main diagonal), $I_{m_k}$ is the $m_k \times
m_k$ identity
matrix,  $V= \left[\begin{array}{cc}
0&v\\v^{*}&0\end{array}\right]$, $v(x)$ is an $m_1 \times m_2$ matrix-valued function (matrix function), and $v^*$ is the conjugate transpose of $v$.
Inverse problem for the closely related Krein system (for the case corresponding to Dirac system
with a continuous scalar potential $v$) was first solved (without proof) in the brilliant note \cite{Kre1}
by M.G. Krein. Some fundamental results for Dirac and related (see Subsection \ref{sub2.2}) canonical systems 
have been obtained in \cite{dB1, dB2, GoKr, SaL2},
see also the references therein.
 For more recent publications and references on Dirac systems see, for instance, 
\cite{AHM, BT, CG, EGNST, FKRS3,  GGHT, GeS, GeZ, MyPuy,  ALS-Sc, SaSaR}.
We note that the inverse problem to recover Dirac system \eqref{1.1} (with a locally bounded potential $v$) from an $m_2\times m_1$
Titchmarsh-Weyl (Weyl) matrix function $\vp(\la)$ was solved in our paper \cite{FKRS3} (see also \cite{SaA3} and \cite[Ch. 2]{SaSaR}). Following a question by F.~Gesztesy, we presented
a solution of the inverse problem for a much more general case, namely, for Dirac systems, the potentials of which have  locally square-integrable entries  \cite{ALS15}.
We speak about  a locally square-integrable matrix function if its entries are locally square-integrable.

{\bf 2.} In the papers \cite{FKRS3, ALS15}, we dealt with $m_2\times m_1$ contractive Weyl matrix functions. For the most important subcase  (of Dirac systems)
$m_1=m_2$,  one may consider Weyl matrix functions belonging to Herglotz class, which opens way to recover
Dirac systems from spectral functions. This inverse problem was treated (in \cite{SaA3}) only  for the case of  locally bounded potentials.

Here, we solve the inverse problem to recover Dirac systems from  spectral functions for the case of   locally square-integrable potentials.
We also characterise the general classes of  distribution matrix functions, which are spectral  functions for some Dirac systems (either on $[0,\ell]$ or on $[0,\infty)$) with locally square-integrable potentials.
It is essential that structured operators, which appear in our paper, are operators
with difference kernels. This allows, in particular, the usage of some important results related to the so called Paley-Wiener sampling ($PW$-sampling) \cite{MP, OrS}.

Further we assume that 
\begin{align} &   \label{1.2}
j = \left[
\begin{array}{cc}
I_{p} & 0 \\ 0 & -I_{p}
\end{array}
\right], \quad V= \left[\begin{array}{cc}
0&v\\v^{*}&0\end{array}\right], \quad v(x)\in \BC^{p\times p}, \quad p\in \BN,
 \end{align} 
where $\BN$ is the set of positive integer numbers, $\BC$ is the complex plane, $\BC^{p\times q}$ is the set of
$p\times q$ complex valued matrices, and the entries of $v(x)$ are locally integrable, which usually suffices for direct spectral problems. 

{\it Notations.} Some notations were introduced above. The notation $\BZ$ stands for the set of integer numbers, $\BR$ stands for the real axis, and $\BC_+$ denotes 
the open upper half-plane of the complex plane $\BC$.  For an interval $\cli$, $|\cli |$ is the length of $\cli$.
By $L^2_p(b,a)$, we denote a Hilbert space
of square-integrable $p\times 1$ vector functions on $(b,a)$ with the scalar product $\langle f,g\rangle_{L^2}=\int_b^{a}g(t)^*f(t)dt$, and
$L^{2}_{p}( d \tau )$ (with some nondecreasing $p\times p$ matrix function $\tau(t)$ on $\BR$)  stands for a Hilbert space
of $p\times 1$ vector functions with the scalar product $\langle f,g\rangle_\tau = \int _{- \infty}^{ \infty}
g(t)^*d \tau (t) f(t)$.  We write $L^2(b,a)$ and  $L^{2}( d \tau )$ instead of $L^2_1(b,a)$ and $L^{2}_{1}( d \tau )$, respectively, in the scalar cases, omitting $p=1$.
The set of bounded operators acting from the Hilbert space $\clh_1$ into the Hilbert space $\clh_2$ is denoted by $B(\clh_1,\clh_2)$
and we write $B(\clh)$ instead of $B(\clh_1,\clh_2)$ in the case $\clh_1=\clh_2=\clh$.  For $S\in B\big(L^2_p(0,a)\big)$, the inequality $S>0$ means that
$\langle S f,f\rangle_{L^2}>0$ if only $f\not=0$.
The set of matrices $\BC^{p\times q}$ was already defined and we denote the column vectors by $\BC^p:=\BC^{p\times 1}$.
Recall that the notations $I_p$ and $v^*$
have been introduced above. The notation $I$ stands for the general type identity operator
and  the Hermitian adjoint of an operator $A$ is denoted by~$A^*$. 
\section{Preliminaries} \label{prel}
\setcounter{equation}{0}
\subsection{Spectral functions of Dirac systems}
{\bf 1.} Let us present some definitions and direct spectral results for  Dirac system \eqref{1.1}, where \eqref{1.2} is always assumed to hold (see \cite[Subsection 1.2.6 and Section 2.1]{SaSaR} for the results below as well as further details
and references).
The operator  $\clh_D$ in $L^2_{2p}(0, \, \infty )$, which corresponds to 
system \eqref{1.1} (so that \eqref{1.1} may be rewritten as $\clh_Df=zf$), is determined by its differential expression $\clh_D$ and domain  of definition ${\cld(\clh_D)}$:
\begin{align}& \label{1.3}
\clh_Df=-\Big(\I  j\frac{d}{dx}+V(x)\Big)f,
\\ &\label{1.4}
{\mathcal{D}} (\clh_D):= \Big\{ f: \,    -\Big(\I j \frac{d}{dx}+V(x)\Big)f \in
L^2_{2p}(0, \, \infty ), \quad [I_p \quad \om]f(0)=0 \Big\}.
\end{align}
It is assumed  that the matrix $\om$
in the initial condition in \eqref{1.4}
is unitary, which is necessary and sufficient for $\clh_D$ being
selfadjoint. Moreover, further in the text we set $\om =-I_p$ since
the results for the general case of unitary $\om$ easily follow (see \cite[pp. 45,46]{SaSaR}).

We denote the fundamental solution of \eqref{1.1} by $u(x,z)$ and normalize it by the condition 
\begin{align}& \label{1.5}
u(0,z)=I_{2p}.
 \end{align}
Put
\begin{equation}         \label{1.6}
U_D= \frac{1}{\sqrt{2} } \int _{0}^{ \infty}\begin{bmatrix}I_{p}
 \quad I_p\end{bmatrix}
 u(x, \overline{ z})^{*} \, \cdot \, dx.
\end{equation}
For functions $f \in  {\cal{D}} (\clh_D) $ with compact support
we easily derive
\begin{align}& \label{1.7}
 (U_D \clh_D f)( z)= z (U_D f)( z ), 
 \end{align}
 that is,
$U_D$ diagonalizes $\clh_D$.
\begin{Dn}     \label{DnSF}
A nondecreasing $p \times p$ matrix function $ \tau $
on the real axis $\BR$
is called a
spectral function of the system
\eqref{1.1} on $[0, \, \infty)$ with the initial condition
\begin{align}& \label{1.8}
\begin{bmatrix}  I_p & -I_p \end{bmatrix}  f(0)=0,
\end{align}
if $U_D$, which is defined by (\ref{1.6}) for $f$
with compact support,
extends to an isometry, also denoted by $U_D$, from
$L^{2}_{2p}(0, \, \infty )$ into
$L^{2}_{p}( d \tau )$. 

Correspondingly, a nondecreasing  matrix function $ \tau $
is called a
spectral function of the system
\eqref{1.1} on the interval $[0, \, \ell]$ $(0<\ell< \infty)$ 
with the initial condition  (\ref{1.8}),
if $U_D$ defined by (\ref{1.6})
restricts
to an isometry from
$L^{2}_{2p}(0, \, \ell)$ into
$L^{2}_{p}( d \tau )$.
\end{Dn}
{\bf 2.} Next, we describe the set of spectral functions. For this, we need more definitions.
\begin{Dn}\label{DefJ} 
A $2p \times  p$ matrix function $\clp(z  )$, which is
meromorphic in $\BC_+$,
is called 
nonsingular with  property-$j$ if the inequalities
\begin{equation}\label{1.9}
\clp(z  )^*\clp(z  )>0,\quad
  \clp(z  )^* j \clp(z  ) \geq 0,  
  \end{equation}
where $j$ is given in \eqref{1.2}, hold in $\BC_+$ $($excluding, possibly, isolated points$)$.
\end{Dn}
Since $\clp$ is meromorphic,  the fact that   the first inequality in \eqref{1.9} holds in one point implies that it holds in all $\BC_+$ (excluding, possibly, isolated points).
We set
\begin{align}  &     \label{1.10}
{\cal{U}}( z )={\cal{U}}_{\ell}( z ):=u(\ell, \overline{ z })^{*}.
\end{align}
If \eqref{1.9} holds, then we have
\begin{align}&      \label{1.11}
\det \Big(\begin{bmatrix}
I_p & I_{p} 
\end{bmatrix}\clu(z)\clp(z)\Big)\not= 0,
\end{align}
and so the M\"obius  transformation
\begin{equation}       \label{1.12}
\vp(z)=\varphi (z , \ell, {\cal P})= \I \begin{bmatrix} I_p & -I_p \end{bmatrix} {\cal U}( z  ){\cal P}(z  ) 
\Big(\begin{bmatrix}I_p & I_p \end{bmatrix}{\cal
U}( z  ){\cal P}( z  )\Big)^{-1}
\end{equation}
is well-defined.
\begin{Nn}\label{NnLFr}
Denote by ${\cal{N}}(\clu)$  the set of M\"obius 
transformations
\eqref{1.12}
where   the matrix functions ${\cal P}$  are nonsingular,
 with 
property-$j$.
\end{Nn}
\begin{Dn} \label{DnW}
M\"obius   transformations $\varphi \in {\cal{N}}(\clu)$  are called Weyl functions of  the Dirac  system
(1.1) on the interval $[0, \ell]$.
\end{Dn}
The $p \times p$ matrix functions $\vp\in \cln(\clu) $ belong to Herglotz class, that is, $\I( \varphi (z)^* - \varphi (z) \geq 0$ in $\BC_+$,
and so they admit a unique Herglotz representation
 \begin{align}  &     \label{1.13}
\varphi (z )= \mu z + \nu +
\int_{- \infty }^{ \infty}
\Big( \frac{1}{t- z } - \frac{t}{1+t^{2}}\Big)d \tau (t),
\end{align}
where $\mu \geq 0$, $\nu= \nu^*$ and $\tau$
is a nondecreasing $p \times p$ distribution matrix function 
such that
 \begin{align}  &     \label{1.13+}
\int_{- \infty }^{ \infty}
(1+t^{2})^{-1}d \tau (t)<\infty .
\end{align}
\begin{Rk}\label{RkSpF}
According to \cite{SaA3} $($see, e.g., Theorems 2.4 and 2.5 in  \cite{SaA3} together with \cite[Proposition 1.49]{SaSaR}$)$,
the set of spectral functions of  Dirac system \eqref{1.1}, \eqref{1.2} on $[0,\ell]$, where the potential $v(x)$ is integrable,
coincides with the set of the distribution matrix functions $\tau$ in the Herglotz representations \eqref{1.13} of  the functions $\vp$ from the set $\cln(\clu)$ of Weyl functions.
\end{Rk}
Note that the matrix functions $\tau(t)$ in the   Herglotz representation of $\vp$ are unique after normalization, that is, after fixing $\tau(t_0)$ for some $t_0\in \BR$ 
and assuming $\tau(t-0)=\tau(t)$. The spectral functions may be normalized in the same way.
\subsection{Canonical systems}\label{sub2.2}
Canonical systems are Hamiltonian systems of the form
\begin{equation}\label{1.14}
dy(x,z ) /dx=\I z JH(x)y(x,z ),\quad
J=\left[\begin{array}{lr}0 & I_p \\ I_p & 0\end{array} \right],
\quad
H(x) \geq 0,
\end{equation}
where the entries of the $2p\times 2p$ Hamiltonian $H(x)$ are locally integrable.

We set
\begin{align}&      \label{1.18}
\Theta:=\frac{1}{\sqrt{2}}\begin{bmatrix}I_p & - I_p\\ I_p & I_p\end{bmatrix} \quad {\mathrm{so\,\, that}}  \quad {\Theta}^*={\Theta}^{-1}, \quad \Theta j \Theta^*=J,
\end{align}
and  partition the fundamental solution $u(\cdot,0)\Theta^*$ of Dirac system (at $z=0$) into two $p\times 2p$ block rows $\b$ and $\g$ and
four $p\times p$ blocks $\b_k$, $\g_k$ $(k=1,2)$:
 \begin{align}  &     \label{1.15}
\b(x)=\begin{bmatrix}\b_1(x) & \b_2(x)\end{bmatrix}=\begin{bmatrix}I_p & 0\end{bmatrix}u(x/2,0)\Theta^*, 
\\    &     \label{1.16}
\g(x)=\begin{bmatrix}\g_1(x) & \g_2(x)\end{bmatrix}=\begin{bmatrix}0 & I_p\end{bmatrix}u(x/2,0)\Theta^*.
\end{align}
Dirac systems are equivalent to the subclass of canonical systems \eqref{1.14}, where 
\begin{align}&      \label{1.17}
H(x)=\b(x)^*\b(x).
\end{align}
{\it Clearly, if Dirac system is considered on the finite interval} $[0,\ell]$, {\it the corresponding
canonical system is considered on} $[0,2\ell]$. Note that $\b$ and $\g$ are introduced as in \cite{SaA3},
which is more convenient for this paper than the slightly different notations in \cite{ALS15}.

According to \eqref{1.1} we have
$$\frac{d}{dx}\Big(u(x/2,\ov{z})^*ju(x/2,z)\Big)=0,$$
which (in view of the normalization \eqref{1.5}) yields 
\begin{align}&      \label{1.20}
u(x/2,\ov{z})^*ju(x/2,z)=u(x/2,z)ju(x/2,\ov{z})^*=j.
\end{align}
Thus, taking also into account \eqref{1.18}--\eqref{1.16}, we obtain  
the following properties of $\b$ and $\g$:
\begin{align}&      \label{1.21}
\b(0)=\frac{1}{\sqrt{2}}\begin{bmatrix}I_p & I_p\end{bmatrix}, \quad \b(x)J\b(x)^*\equiv I_p, \quad \b(x)J\g(x)^*\equiv 0,
\\ &      \label{1.22}
\g(0)=\frac{1}{\sqrt{2}}\begin{bmatrix}-I_p & I_p\end{bmatrix}, \quad \g(x)J\g(x)^*\equiv -I_p.
\end{align}
Moreover, relations \eqref{1.1}, \eqref{1.15} and \eqref{1.16} imply that
\begin{align}&      \label{1.23}
\b^{\prime}(x)=\frac{\I}{2}v(x/2)\g(x)\quad \g^{\prime}(x)=-\frac{\I}{2}v(x/2)^*\b(x)  \quad \Big(\b^{\prime}=\frac{d}{dx} \b\Big),
\end{align}
and so (in view of the last relations in \eqref{1.21} and \eqref{1.22}) we have
\begin{align}&      \label{1.23+}
\b^{\prime}(x)J \b(x)^*\equiv 0, \quad \g^{\prime}(x)J \g(x)^* \equiv 0,
\\ &      \label{1.24}
v(x/2)=2\I \b^{\prime}(x)J\g(x)^*.
\end{align}
For the fundamental solution $W$ of the system \eqref{1.14}, \eqref{1.15}, \eqref{1.17} normalized by
$W(0,z)=I_{2p}$, we easily obtain (see also \cite[(1.4) and (4.5)]{SaA3}) that
\begin{align}&      \label{1.19}
W(x,z)=\exp\{\I z x/2\}           \Theta \, u(x/2,0)^{-1}u(x/2,z)\Theta^*.
\end{align}
\subsection{Linear similarity result}
The similarity result below (a particular case of \cite[Proposition 2.1]{ALS15}) plays an essential role
in solving the inverse problem. 
Namely, we consider the conditions of similarity of a Volterra operator $\clk$ to the integration operator~$\cla$:
\begin{equation} \label{1.25}
\clk:=F(x)\int_0^x G(t)\, \cdot \,dt, \quad \cla:=\int_0^x \, \cdot \,dt \quad \big(\clk,\cla\in B\big(L^2_p(0, \,\eT)\big)\big).
\end{equation}
where $F$ and $G$ are differentiable $p \times 2p$ and $2p \times p$, respectively, matrix functions.
\begin{Pn} \label{PnSim} Let $F$ and $G$ be differentiable and satisfy the identity
\begin{align}\label{1.26}
F(x)G(x) \equiv I_p, \quad 0 \leq x \leq \eT,
\end{align}
and assume that the entries of $F^{\prime}$ and $G^{\prime}$ belong to $L^2(0, \, \eT)$.

Then, the operator $\clk$ defined by \eqref{1.25}
is linearly similar to the operator of integration $\cla$. More precisely, $\clk=E \cla
E^{-1}$ where $E \in B\big(L^2_p(0, \, \eT)\big)$ is a lower triangular operator of the form
\begin{equation}
\label{1.27} E= \rho(x)\left(I+\int_0^x N(x,t)\, \cdot \,dt\right), \quad \frac{d}{dx}\rho=F'G\rho, \quad \rho(0)=I_p,
\end{equation}
and the matrix functions $\rho$, $\rho^{-1}$ and $N$ are measurable and uniformly bounded.
Moreover, the operators $E^{\pm 1}$ map differentiable functions with a square-integrable derivative on $[0,\eT]$
into differentiable functions with a square-integrable derivative on $[0,\eT]$.
\end{Pn}

Let us assume that the entries of $v(x)$ are locally square-integrable, which means that they are square-integrable on $[0,\ell]$. We  set $\eT=2\ell$
and will be interested in the operator $\clk=\b(x)J \int_0^x \b(t)^*\, \cdot \,dt\in  B\big(L^2_p(0, \, 2\ell)\big)$. Clearly, in view of \eqref{1.21} and \eqref{1.23+}, $F(x)=\b(x)$
and $G(x)=J\b(x)^*$ satisfy the requirements of  Proposition \ref{PnSim} and $\rho(x)\equiv I_p$ in \eqref{1.27}. 
Thus, $\clk$ is linearly similar to the integration. Note that the relation $\b J\b^*\equiv I_p$ required in \eqref{1.26} 
 is immediate from $\b(0)J\b(0)^*=I_p$ and $\b^{\prime}J\b^*\equiv 0$.
\begin{Cy}\label{CySim} Let  $\b(x)$ be a differentiable $p\times 2p$ matrix function on $[0,2\ell]$ with the square-integrable derivative.
Assume that $\b(0)J\b(0)^*=I_p$ and $\b^{\prime}J\b^*\equiv 0$.
$($In particular, the conditions of the corollary are satisfied for $\b$ generated via \eqref{1.15} by a Dirac system with the square-integrable on $[0,\ell]$
potential $v.)$
Then, we have
\begin{equation}
\label{1.28} K:=\I \b(x)J \int_0^x \b(t)^*\, \cdot \,dt=\wt E A \wt E^{-1},  \end{equation}
where $A,\wt E\in B\big(L^2_p(0, \, 2\ell)\big)$,
\begin{equation}
\label{1.29} A:=\I \int_0^x \, \cdot \,dt, \quad  \wt E= I+\int_0^x \wt N(x,t)\, \cdot \,dt, 
\end{equation}
and the similarity transformation operator $\wt E$ has all the properties described in Proposition \ref{PnSim}.
\end{Cy}
\section{Inverse problem} \label{inv}
\setcounter{equation}{0}
{\bf 1.}  Let us assume that Dirac system \eqref{1.1} with a   square-integrable potential is given.
Then (as stated in Corollary \ref{CySim}), there is a similarity transformation operator  such that
it has all the properties described in Proposition~\ref{PnSim} and that \eqref{1.28} holds.
This operator is not uniquely defined. One may, for instance, fix first such an  operator constructed
in \cite[Proposition 2.1]{ALS15} and denote it by $\wt E$. Then, we consider the operator
\begin{equation}\label{3.1} 
E=\wt E E_0, \quad E_0=\frac{1}{\sqrt{2}}I + \int_0^x \cle_0(x-t) \cdot dt, \quad 
\cle_0(x):=\big(\wt E^{-1}\b_2\big)^{\prime}(x),
\end{equation}
where the matrix function $\b_2$ is introduced in \eqref{1.15} and $\wt E^{-1}$ is applied  to the columns of $\b_2$.
Recall that $A,K$ and $E$ are acting in $L^2_p(0,2\ell)$.
\begin{La} \label{LaSim} Let Dirac system  \eqref{1.1}  with    a square-integrable potential be given on some interval $[0,\ell]$,
let $\b(x)$ have the form \eqref{1.15} and $E$ have the form \eqref{3.1}. Then, 
\begin{equation}
\label{3.2} K:=\I \b(x)J \int_0^x \b(t)^*\, \cdot \,dt= E A  E^{-1},  \end{equation}
where
\begin{equation}
\label{3.3} A:=\I \int_0^x \, \cdot \,dt, \quad   E=\frac{1}{\sqrt{2}} I+\int_0^x  N(x,t)\, \cdot \,dt, 
\end{equation}
and $N(x,t)$ is a Hilbert-Schmidt kernel.
Moreover,  the operators $E^{\pm 1}$ map differentiable functions with a square-integrable derivative
into differentiable functions with a square-integrable derivative.
Finally, we have the equality
\begin{equation}
\label{3.4} E^{-1}\b_2(x)\equiv I_p.
\end{equation}
\end{La}
\begin{proof}  The second relation in \eqref{1.29} (where  $\wt N$ is bounded) and formulas \eqref{3.1} imply the representation \eqref{3.3} of $E$,
 where $N(x,t)$ is a Hilbert-Schmidt kernel. The second and third relations in \eqref{3.1} yield also
 \begin{equation}\label{3.5} 
 E_0 I_p=\big(\wt E^{-1}\b_2\big)(x)-\big(\wt E^{-1}\b_2\big)(0)+\frac{1}{\sqrt{2}}I_p.
\end{equation}
Moreover,  it follows from \eqref{1.21} and \eqref{1.29} that $\big(\wt E^{-1}\b_2\big)(0)=\b_2(0)=\frac{1}{\sqrt{2}}I_p$.
Hence, we rewrite \eqref{3.5} as $E_0 I_p=\big(\wt E^{-1}\b_2\big)(x)$. Now, \eqref{3.4} follows from the definition of $E$
(i.e., from the first equality in \eqref{3.1}).

It is easy to see that  $E_0$ commutes with $A$.  Indeed,  consider the relation
\begin{align}&\nn
\int_0^x\int_0^t g(t-\xi)f(\xi)d\xi dt =\int_0^x\int_{\xi}^x g(t-\xi)dt f(\xi)d\xi
\\ &\label{C}
=\int_0^x\int_{0}^{x-\xi} g(u)du \,\, d \int_0^{\xi}f(\zeta)d\zeta
=\int_0^x g(x-t)\int_{0}^t f(\xi)d\xi dt,
\end{align}
where $f\in L^2_p(0,\ell)$, $g$ is some $p\times p$ matrix function with columns  belonging
to $L^2_p(0,\ell)$ and the last equality is obtained using integration by parts. Clearly,
\eqref{C} yields $AE_0=E_0A$. Hence, similarity formula \eqref{3.2} follows from \eqref{1.28} and the first
equality in \eqref{3.1}.

Now, taking into account the properties of $\wt E$, it suffices to show that
 the operators $E_0^{\pm 1}$ map differentiable functions with a square-integrable derivative
into differentiable functions with a square-integrable derivative in order to prove the lemma.

The mentioned above property of the mapping $E_0$ is immediate from the fact that $E_0A=AE_0$ and that $E_0 I_p=\big(\wt E^{-1}\b_2\big)(x)$
consists of differentiable vector functions (columns) with the square-integrable derivatives. The corresponding property of $E_0^{-1}$
follows from the rewriting of $E_0A=AE_0$ and  $E_0 I_p=\big(\wt E^{-1}\b_2\big)(x)$ in the forms
 \begin{align} \nn
E_0^{-1}A=AE_0^{-1}, \quad \frac{1}{\sqrt{2}} E_0^{-1}I_p&=I_p+\I E_0^{-1}A\big(\wt E^{-1}\b_2\big)^{\prime}
\\ & \label{z0}
=I_p+\I A E_0^{-1}\big(\wt E^{-1}\b_2\big)^{\prime}.
\end{align}
\end{proof}
{\bf 2.} Next, we introduce the operator identity 
\begin{equation}\label{3.6} 
AS-SA^*=\I \Pi J \Pi^* \quad (S\in B\big(L^2_p(0,2\ell)\big),
\end{equation}
and the corresponding structured operator $S$, both playing an essential role in solving inverse problem.
Namely, for $K$ given in \eqref{3.2}, it is easy to see that $K-K^*=\I \b(x)J\int_0^{2\ell}\b(t)^* \, \cdot \, dt$, which we rewrite in the form
\begin{equation}\label{3.7} 
E A  E^{-1}-\big(E^{-1}\big)^* A^*  E^*=\I \b(x)J\int_0^{2\ell}\b(t)^* \, \cdot \, dt.
\end{equation}
We multiply both parts of \eqref{3.7} by $E^{-1}$ from the left and  by $\big(E^{-1}\big)^*$ from the right.
Thus, we obtain the identity \eqref{3.6}, where
\begin{align}& \label{3.8} 
S=  E^{-1}\big(E^{-1}\big)^* , \quad \Pi=\begin{bmatrix}\Phi_1 & \Phi_2\end{bmatrix}, \quad \Phi_1, \Phi_2 \in B\big(\BC^{p}, L^2_p(0, \, 2\ell)\big),
\\ & \label{3.9} 
\Phi_kg=\Phi_k(x)g, \quad \Phi_k(x):=\big(E^{-1}\b_k\big)(x) \quad (g\in \BC^p, \,\, k=1,2).
\end{align}
Clearly, $S$ is a bounded and strictly positive operator.
In view of  \eqref{3.4} and \eqref{3.9}, the operator $\Phi_2$ is an embedding of $\BC^p$ into $L^2_p(0, \, 2\ell)$, that is,
\begin{equation}\label{3.10} 
\Phi_2(x)\equiv I_p.
\end{equation}
\begin{Rk}\label{RkPhi1} Further in the text, we will need the equality
\begin{equation}\label{3.10+} 
\Phi_1(0)= I_p,
\end{equation}
which easily follows from Lemma \ref{LaSim} with its proof and from \eqref{3.9}.
\end{Rk}
We introduce the transfer matrix function in Lev Sakhnovich form~\cite{SaL1, SaL2, SaL21}
\begin{align} &      \label{3.11}
w_A(z):=I_{2p}+\I z J\Pi^*S^{-1}(I-zA)^{-1}\Pi.
\end{align} 
We will also use  the reductions of the operators $A$ and $S$ above:
\begin{align} &      \label{3.12}
\big(P_{\xi}f\big)(x)=f(x) \quad (0<x<\xi\leq 2\ell), \quad P_{\xi}\in B\Big(L^2_{p}(0, \, 2\ell), \, L^2_{p}(0, \, \xi)\Big),
\\ &      \label{3.13}
 A_{\xi}:=P_{\xi}AP_{\xi}^*, \quad S_{\xi}:=P_{\xi}SP_{\xi}^*,
\end{align}
where $P_{\xi}^*$ is a natural embedding of  $L^2_{p}(0, \, \xi)$ into $L^2_{p}(0, \, 2\ell)$.
Since $A$ is a lower triangular operator and the operator identity \eqref{3.6} holds, we have
\begin{align} &      \label{3.14-}
P_{\xi}A=P_{\xi}AP_{\xi}^*P_{\xi}=A_{\xi}P_{\xi}, \quad A_{\xi}S_{\xi}-S_{\xi}A_{\xi}^*=\I P_{\xi}\Pi J \Pi^*P_{\xi}^*.
\end{align} 
The matrix function $w_A(\xi,z)$ corresponding to the reductions \eqref{3.13} has the form
\begin{align} &      \label{3.14}
w_A(\xi,z):=I_{2p}+\I z J\Pi^*P_{\xi}^*S_{\xi}^{-1}(I-zA_{\xi})^{-1}P_{\xi}\Pi, \quad 0<{\xi}\leq 2\ell .
\end{align}
Similar to \eqref{3.14-}, one can show (see \cite[p. 560]{ALS15}) that
\begin{align} &      \label{3.15-}
S_{\xi}^{-1}=E_{\xi}^*E_{\xi}, \quad E_{\xi}:=P_{\xi}EP_{\xi}^*.\end{align} 
Relations \eqref{3.9} yield $E \Pi g=\b(x) g$ $(g\in \BC^{2p})$ on $[0,2\ell]$ and so 
\begin{align} &      \label{3.14+}
E_{\xi}P_{\xi} \Pi g=\b(x) g \quad (0<x<\xi\leq 2\ell).
\end{align}
Thus, taking into account \eqref{3.8}, \eqref{3.9} and \eqref{3.15-}, we derive
\begin{align} &      \label{3.15}
\frac{d}{d\xi}\Big(\Pi^*P_{\xi}^*S_{\xi}^{-1}P_{\xi}\Pi\Big)=\b(\xi)^*\b(\xi)=H(\xi).
\end{align} 
Thus, according to the Continuous Factorization Theorem \cite[p. 40]{SaL2} or its corollary \cite[Theorem 1.20]{SaSaR}
we obtain a transfer matrix function representation of the fundamental solution $W(x,z)$ of the canonical system \eqref{1.14}, \eqref{1.17}:
\begin{align} &      \label{3.16}
W(\xi,z)=w_A(\xi,z).
\end{align} 
Weyl and spectral functions of the canonical systems generated by the so called $S$-nodes (i.e., by the triples $\{A, S,\Pi\}$)  are expressed in terms of the M\"obius transformations
\begin{align} &      \label{3.17}
\phi(z)= \I \begin{bmatrix} I_p & 0\end{bmatrix} {\A}( z  )\wh {\clp}(z  ) 
\Big(\begin{bmatrix}0 & I_p \end{bmatrix}{
\A}( z  )\wh {\clp}( z  )\Big)^{-1}, \end{align} 
where $\wh {\cal P}(z  )$ are $2p \times p$ nonsingular matrix functions with property-$J$ and 
\begin{align} &      \label{3.18}
 {\A}( z  )=\A_{\ell}(z):=w_A(2\ell,\ov{z})^*.
\end{align} 
The family of $p\times p$ matrix functions $\phi(z)$ (see \eqref{3.17}) is denoted  by $\cln(\A)$.
Relations \eqref{1.10}, \eqref{1.19}, \eqref{3.16} and \eqref{3.18} yield
\begin{align} &      \label{3.19}
 {\A}( z  )=\exp\{-\I z \ell\}\Theta \, \clu(z)\big(u(\ell,0)^{-1}\big)^*\Theta^*.
\end{align} 
In view of \eqref{1.18} and \eqref{1.20}, both transformations $\clp(z)=\Theta^*\wh{\clp}(z)$ and $\clp(z)=\big(u(\ell,0)^{-1}\big)^*\Theta^*\wh{\clp}(z)$
map the class of nonsingular matrix functions
with property-$J$ onto the class of nonsingular matrix functions
with property-$j$. Thus, comparing M\"obius transformations \eqref{1.12} and \eqref{3.17} and taking into account \eqref{3.19}
and the definition of $\Theta$ in \eqref{1.18}, we derive
\begin{align} &      \label{3.20}
\cln(\clu)=\cln(\A).
\end{align} 
According to \cite[Ch. 4]{SaL2} (see also \cite[Theorem 1.27]{SaSaR}), we have $\mu=0$ in the Herglotz representations \eqref{1.13}
of  the functions $\phi \in \cln(\A)$ or, equivalently, of the functions $\vp\in  \cln(\clu)$. Moreover, the corresponding spectral functions $\tau$
and matrices $\nu$ provide representations of the operators $S$ and $\Phi_1$ introduced in \eqref{3.8} and \eqref{3.9}:
\begin{align}&\label{3.21}
S=  \int_{- \infty}^{\infty}(I-t A )^{-1} \Phi_{2} ( d \tau(t))
\Phi_{2}^* (I-t A^*)^{-1},
\\ &\label{3.22}
\Phi_1= \I \left( \Phi_{2} \nu- \int_{- \infty}^{\infty} \left(
A(I-t A )^{-1} + \frac{t}{t^2+1}I \right) \Phi_{2} d \tau(t)
\right),
\end{align}
where the right-hand sides of  \eqref{3.21} and \eqref{3.22} weakly converge.
\begin{Rk} \label{RkRes} Recall that $\Phi_2(x)\equiv I_p$  in our case and note the resolvent of the operator $A$ above has  a simple explicit form:
\begin{align}\label{3.23} &
(I-zA)^{-1}=I+\I z\int_0^x\E^{\I z(x-\xi)}\, \cdot \, d\xi, \quad (I-zA)^{-1}\Phi_2 g= \E^{\I z x}g \quad (g\in \BC^p).
\end{align}
\end{Rk}
\begin{Rk}\label{RkConv} The weak convergence of the right-hand side of \eqref{3.22} follows from the weak convergence
of the right-hand side of \eqref{3.21} under quite general conditions  \cite[p. 56]{SaL2} and the proof is especially simple in our case.

Indeed, it is easy to see that the equality
\begin{align} &      \label{z7}
A(I-tA)^{-1}+\frac{t}{1+t^2}I=(1+t^2)^{-1}(A+tI)(I-tA)^{-1}
\end{align} 
is valid. In view of the second equality in \eqref{3.23}, the weak convergence in \eqref{3.21} means that 
\begin{align}\label{z8} &
\Phi_{2}^* (I-t A^*)^{-1}f=\int_0^{\ell}\E^{-\I t x}f(x)dx\in L^2_p(d\tau) \quad {\mathrm{for\,\, any}} \quad f\in L^2_p(0,\ell).
\end{align}
Thus, taking into account \eqref{1.13+}, we see that  the integrals of the form
\begin{align} &      \nn
\int_{-\infty}^{\infty}\int_0^{\ell}f(x)^*\E^{\I t x}dx\big(d\tau(t)\big)t(1+t^2)^{-1}g, \quad g\in \BC^p
\end{align}
converge, that is, weakly converges the integral
\begin{align} &      \label{z9}
\int_{-\infty}^{\infty} t(1+t^2)^{-1}(I-tA)^{-1}\Phi_2 d\tau(t).
\end{align}
From \eqref{1.13+}, it is also easy to see that the integral
\begin{align} &      \label{z10}
\int_{-\infty}^{\infty} (1+t^2)^{-1}A(I-tA)^{-1}\Phi_2 d\tau(t)=\int_{-\infty}^{\infty} \big(t(1+t^2)\big)^{-1}\big(\E^{\I t x}-1\big)I_p d\tau(t)
\end{align}
converges $($strongly$)$ as well. Formula \eqref{z7} and convergences in \eqref{z9} and \eqref{z10} imply that the right-hand side
of \eqref{3.22} weakly converges.
\end{Rk}
 We express $\nu$ in terms of the spectral function $\tau$, that is,
in terms of the integral part on the right-hand side of \eqref{3.22}. The corresponding operator is denoted by $\Lam$,
where $\Lam g=\Lam(x)g$ $(\Lam(x)\in \BC^{p\times p}, \,\, g\in \BC^p)$. More precisely, in view of \eqref{3.10}, \eqref{3.10+} and \eqref{3.22}
we have
\begin{align}\label{3.24} &
\nu=\Lam(0)-\I I_p, \,\, {\mathrm{where}}\,\, \Lam(x)=\int_{- \infty}^{\infty} \left(
A(I-t A )^{-1} + \frac{t}{t^2+1}I \right) I_p\, d \tau(t).
\end{align}
{\bf 3.} Now, we formulate our main theorem on inverse problem, which is already proved above. 
Note that the recovery of $\b(x)$ and $\g(x)$ from $S$ and $\Pi$
is described in the appendix.
\begin{Tm} \label{TmInvPr} Let Dirac system \eqref{1.1}, \eqref{1.2} be given on the interval $[0,\ell]$
and have a square-integrable $($on $[0,\ell])$  $p\times p$ matrix potential $v(x)$.
Then, this system may be uniquely recovered from its spectral  function $\tau$ in three steps.
First, we recover operators $S$ and $\Pi=\begin{bmatrix}\Phi_1 & \Phi_2\end{bmatrix}$ using relations \eqref{3.21}, \eqref{3.22}, \eqref{3.24} and the
equality $\Phi_2g\equiv g$ $(g\in \BC^p)$.  Next, we recover $\b(x)$ and $\g(x)$ using relations \eqref{3.15}, \eqref{A2}--\eqref{A4}
and \eqref{A6},\eqref{A7}. Finally, the potential $v$ is given by the formula \eqref{1.24}:
$$v(x/2)=2\I \b^{\prime}(x)J\g(x)^*.$$
\end{Tm}
\begin{Rk} A way to directly recover $\g$ from $S$ and $\Pi$ $($without using $\b$ as in \eqref{A6}$)$ is given by formula \eqref{A9}.
After that, one may recover $\b$ from $\g$   similar to the recovery of $\g$ from $\b$ (see Proposition \ref{PnBG}).
One may also use \eqref{3.15} and \eqref{A2}--\eqref{A4} to recover $\b$ independently. Thus, we have several ways to recover
$\b$ and $\g$.
\end{Rk}
{\bf 4.}
Let us consider Weyl and spectral functions on the semiaxis $[0,\infty)$. For this purpose, we need different intervals $[0,\ell]$ ($\ell$ is not fixed
anymore) and we write $\cln(\clu_{\ell})$ and $\cln(\A_{\ell})$ instead of  $\cln(\clu)$ and $\cln(\A)$. (See  \eqref{1.10} and \eqref{3.18} for $\clu_{\ell}$
and $\A_{\ell}$.) Together with Hamiltonians $H$ of the form \eqref{1.17}, canonical systems  \eqref{1.14} with Hamiltonians
\begin{align}\label{3.25} &
\wh H(x)=2H(2x)-J=\Th u(x,0)^*u(x,0)\Th^*
\end{align}
have been studied in \cite{SaA3}. The fundamental solutions of these systems will be denoted by $\wh W(x,z)$ (and we normalize  them by $\wh W(0,z)=I_{2p}$).
The fundamental solutions of the canonical systems with Hamiltonians \eqref{1.17} and \eqref{3.25} are connected by the relation
\begin{align}\label{3.26} &
W(x,z)=\exp\{\I z x/2\}\wh W(x/2,z).
\end{align}
In view of \eqref{3.20} and \eqref{3.26}, for $\wh \A_{\ell}(z)=\wh W(\ell,\ov{z})^*$ we have
\begin{align} &      \label{3.27}
\cln(\clu_{\ell})=\cln(\A_{\ell})=\cln(\wh \A_{\ell}),
\end{align} 
where the family of  functions $\phi \in \cln(\wh \A_{\ell})$ is obtained via the substitution of $\wh \A_{\ell}(z)$ instead of $\A(z)=\A_{\ell}(z)$ into the right-hand side of
\eqref{3.17}.
The families in \eqref{3.27} are decreasing (i.e., $\cln(\clu_{\ell}) \subseteq  \cln(\clu_r)$ for $\ell>r$).
Moreover, according to \cite[Corollary 2.12]{SaA3} (and to its proof)  the intersection $ \bigcap_{\ell<\infty}\cln(\wh \A_{\ell})$
consists of one and only one function. (Note that the set of matrix functions $\pi \phi(z)$ where $\phi\in \cln(\wh \A_{\ell})$ is denoted
by $\cln(H,\ell)$ in \cite{SaA3}.) The function, which belongs to the intersection $\bigcap_{\ell<\infty}\cln(\wh \A_{\ell})$ or, equivalently,
$\bigcap_{\ell<\infty}\cln(\clu_{\ell})$ we denote by $\vp_{\infty}$. 
\begin{Rk}\label{RkInfty}
Thus, by virtue of  \cite[Corollary 2.12]{SaA3}, we have 
\begin{align} &      \label{3.28}
\vp_{\infty}=\displaystyle{\bigcap_{\ell<\infty}\cln(\clu_{\ell})}
\end{align} 
for a Dirac system with the locally integrable on $[0,\infty)$ potential $v(x)$. In view of Remark \ref{RkSpF} and \eqref{3.28}, there is a unique 
spectral function of the Dirac system with the locally integrable on $[0,\infty)$ potential and this spectral function is the distribution function
$($denoted by $\tau_{\infty})$ in the Herglotz representation of the matrix function $\vp_{\infty}$. Formula \eqref{3.28} may serve as a definition
of the Weyl function on the semiaxis. According to Definition 5.1 and Proposition~5.2 in \cite{SaA3}, the requirement below is an equivalent definition
of the Weyl function:
\begin{equation} \label{3.29}
\int_0^\infty \left[ \begin{array}{lr} I_p &  \I \varphi (z )^*
\end{array} \right]
  \Theta  u(x, z )^*
 u(x, z )\Theta ^*
 \left[ \begin{array}{c}
I_p \\ - \I \varphi (z ) \end{array} \right] dx < \infty , \quad z 
\in {\BC}_+.
\end{equation}
$($Recall that $\Th$ has the form \eqref{1.18}.$)$
\end{Rk}
\begin{Cy}\label{CyInv} Dirac system \eqref{1.1}, \eqref{1.2} with a locally square-integrable on $[0,\infty)$  matrix potential $v(x)$
is uniquely recovered from its spectral  function $\tau_{\infty}$ using the procedure described in Theorem  \ref{TmInvPr} for the
finite intervals.
\end{Cy}
Using the first two equalities in \eqref{3.14-}, we obtain
\begin{align} &      \label{3.30}
 P_{\xi}(I-tA)^{-1}=(I-tA_{\xi})^{-1}P_{\xi}.
\end{align} 
\begin{Rk}\label{Rklr}
In view of \eqref{3.30}, we easily derive from \eqref{3.10}, \eqref{3.21}, \eqref{3.22} and \eqref{3.24} that $S_{\xi}$ and $P_{\xi}\Pi$
do not depend on the choice of the interval $[0,\ell]$, on which we consider Dirac system $($and, correspondingly, of the interval $[0,2\ell]$ for the operators $S$ and $\Pi)$
as long as $\xi\leq 2\ell$. In particular, it follows that $\b(x)$ and $\g(x)$ recovered in Theorem  \ref{TmInvPr} and Corollary \ref{CyInv} do not depend on the choice of $\ell$ $(\ell \geq x/2)$.
Thus, the recovery of Dirac system on the semiaxis in Corollary \ref{CyInv} is well defined, unique and depends only on $\tau$.
\end{Rk}
\section{On the sets of spectral functions} \label{set}
\setcounter{equation}{0}
{\bf 1.} According to the considerations of Section \ref{inv} and to Lemma \ref{LaSim}, formulas \eqref{1.13+}, \eqref{3.8}, \eqref{3.9},  Remark \ref{RkConv} and Theorem \ref{TmInvPr}, in particular,
each spectral function $\tau(t)$ of  a Dirac system with a square-integrable $p\times p$ matrix potential $v(x)$ on $[0,\ell]$ has the following properties (or, equivalently, satisfies the following
conditions). \\
{\bf Properties $(i)-(iii)$}:

$(i)$ The matrix function $\tau$ is  nondecreasing and the inequality \\ $\int_{- \infty }^{ \infty}
(1+t^{2})^{-1}d \tau (t)<\infty $ is valid.

$(ii)$ The right-hand side of the equality
\begin{align}&\label{C1}
S=  \int_{- \infty}^{\infty}(I-t A )^{-1} \Phi_{2} ( d \tau(t))
\Phi_{2}^* (I-t A^*)^{-1},
\end{align}
where $A$ is given in \eqref{3.3}, $\Phi_2\in B\big(\BC^p,L^2_p(0,2\ell)\big)$ and $\Phi_2g\equiv g$, weakly converges and the corresponding operator $S$ is bounded in $L^2_p(0,2\ell)$
and positive (i.e., $S>0$). Moreover, there is an inverse operator  $S^{-1}\in B\big(L^2_p(0,2\ell)\big)$.

Before introducing properties $(iii)$ we need some preparations.
Recall (see Remark \ref{RkConv}) that the convergences of the integrals in $(i)$ and $(ii)$ yield a weak convergence of the integral
\begin{align}&\label{C2-}
\Lam=  \int_{- \infty}^{\infty} \left(
A(I-t A )^{-1} + \frac{t}{t^2+1}I \right) \Phi_{2} d \tau(t), \,\, {\mathrm{where}}\,\, \Lam g=\Lam(x)g 
\end{align}
for $g\in \BC^p$. Since $\Lam(x)$ is continuous (and even differentiable together with $\Phi_1(x)$), the matrix $\nu$  given by \eqref{3.24} is well defined:
\begin{align}&\label{C2n}
\nu=\Lam(0)-\I I_p.
\end{align}
Note that instead of requiring the continuity of  $\Lam(x)$ one may require that $\Phi_1(x)$ is continuous (for each $\nu$ or, equivalently, for some $\nu$), where
\begin{align}&\label{C2}
\Phi_1= \I \big( \Phi_{2} \nu- \Lam\big), \quad \Phi_1 g= \Phi_1(x)g \quad (g\in \BC^p).
\end{align}
 The next properties of $\tau$ are expressed in terms of $\Phi_1(x)$ and $\nu$.

$(iii)$ The entries of the the $p\times p$ matrix function $\Phi_1(x)$ are differentiable and the columns of $\Phi_1^{\prime}(x)$ belong to $L^2_p(0,2\ell)$.  
Moreover, in view of   \eqref{1.13} $\nu$ appearing in \eqref{C2n} and \eqref{C2}   is self-adjoint:
\begin{align} &      \label{C2++}
\nu=\nu^*.
\end{align} 
\begin{Rk}\label{RkBST} The weak convergence of the integral on the right-hand side of \eqref{C1} implies that the operator $S$
in \eqref{C1} is well defined and bounded. The result follows from Banach-Steinhaus theorem $($i.e., from the uniform boundedness
principle$)$.
\end{Rk}
Taking into account \eqref{C2n} and \eqref{C2} (see also \eqref{3.10+}), we have:
\begin{align} &      \label{C2+}
\Phi_1(0)= I_p.
\end{align} 
{\bf 2.} We will show that the properties $(i)-(iii)$ are not only necessary but sufficient for $\tau$ to be a spectral function of  a Dirac system with a square-integrable 
potential.
\begin{Tm}\label{TmSets} The properties $(i)-(iii)$ presented
above are necessary and sufficient for a given $p\times p$ matrix function $\tau(t)$ on $\BR$ to be a spectral function of  some Dirac system on $[0,\ell]$ with a square-integrable 
potential.
\end{Tm}
In order to prove Theorem \ref{TmSets}, we need several auxiliary statements.
\begin{La}\label{LaId} Let a given  $p\times p$ matrix function $\tau(t)$ on $\BR$ have the properties $(i)-(iii)$.
Then, the operators $S$ and $\Phi_1$ determined by \eqref{C1} and \eqref{C2} satisfy the operator identity
\begin{align}\label{C3} &
AS-SA^*=\I (\Phi_1\Phi_2^*+\Phi_2\Phi_1^*)=\I \Pi J\Pi^*, \quad \Pi:=\begin{bmatrix}\Phi_1 &\Phi_2\end{bmatrix}.
\end{align}
Moreover, $S$ is a so called operator with difference kernel $($or convolution operator$):$
\begin{align} \label{C10} &
S=\frac{d}{dx}\int_0^{2\ell}s(x-t) \, \cdot \, dt=2I+\int_0^{2\ell}s^{\prime}(x-t) \, \cdot \, dt,
\end{align}
where
\begin{align} \label{C9} &
s(x)=\Phi_{1}(x) \,\, {\mathrm{for}}\,\, x>0, \quad s(x)=-s(-x)^* \,\, {\mathrm{for}}\,\, x<0. 
\end{align}
\end{La}
\begin{proof} 
It is stated in \cite[p. 3]{SaL1+} that representations \eqref{C1} and \eqref{C2}, \eqref{C2-} yield the operator identity \eqref{C3} (as in \eqref{3.6}).
Although the result is clear enough, it could be convenient to prove it in greater detail.
Indeed, let us set
\begin{align}\label{C4} &
\wh S=  \int_{\G}(I-t A )^{-1} \Phi_{2} ( d \tau(t))
\Phi_{2}^* (I-t A^*)^{-1}, 
\\ \nn &
 \Psi_1=\int_{\G}(1/t)(I-t A )^{-1} \Phi_{2} ( d \tau(t))
\Phi_{2}^* (I-t A^*)^{-1}, \quad 
\Psi_2=\int_{\G}(1/t)\Phi_{2} ( d \tau(t))
\Phi_{2}^*,
\\ & \label{C5}
 \wh \Phi_1= -\I  \int_{\G} \left(
A(I-t A )^{-1} + \frac{t}{t^2+1}I \right) \Phi_{2} d \tau(t),
\end{align}
where $\G=(-T_1,\ve)\cup(\ve,T_2) \quad (0<\ve<T_k, \,\, k=1,2)$.
Then, we have
\begin{align}\nn
A\wh S &=\int_{\G}\big(A-(1/t)I\big)(I-t A )^{-1}\Phi_{2} ( d \tau(t))
\Phi_{2}^* (I-t A^*)^{-1}+\Psi_1
\\ \nn &
=-\int_{\G}(1/t)\Phi_{2} ( d \tau(t))
\Phi_{2}^* (I-t A^*)^{-1} +\Psi_1
\\ \label{C6} &
=-\Psi_2-\int_{\G}\Phi_{2} ( d \tau(t))
\Phi_{2}^* (I-t A^*)^{-1}A^*+ \Psi_1.
\end{align}
In a similar way, we derive
\begin{align} \label{C7} &
\wh SA^* =\Psi_1-\Psi_2-\int_{\G}A(I-t A )^{-1}\Phi_{2} ( d \tau(t))
\Phi_{2}^*.
\end{align}
From \eqref{C5}--\eqref{C7} we obtain
\begin{align}  \nn
A\wh S-\wh SA^* &=\int_{\G}A(I-t A )^{-1}\Phi_{2} ( d \tau(t))
\Phi_{2}^*-\int_{\G}\Phi_{2} ( d \tau(t))
\Phi_{2}^* (I-t A^*)^{-1}A^*
\\ &  \label{C8}
=\I (\wh \Phi_1\Phi_2^*+\Phi_2\wh\Phi_1^*).
\end{align}
After applying some limits, formula \eqref{C8} yields the identity \eqref{C3}.

Now, it is immediate from the properties $(iii)$ and from \cite[Corollary~D.4]{SaSaR} that the unique operator $S$ satisfying
identity \eqref{C3} (with the considered here $A$ and $\Phi_2$ and fixed $\Phi_1$)  has the form \eqref{C10}, \eqref{C9}.
In particular, the second equality in \eqref{C10} follows from  \eqref{C2+}.
\end{proof}
{\bf 3.}
Using some results from \cite{ALS-HLB} we will prove the proposition below on the factorisation of $S$.
\begin{Pn}\label{Pn} Let a convolution operator $S> 0$ of the form \eqref{C10}, \eqref{C9} have a bounded
inverse $\big($i.e., let $S^{-1}\in B\big(L^2_p(0,2\ell)\big)\big)$ and  assume that   $\Phi_1(x)$ in \eqref{C9}
has the properties $(iii)$. Then, $S^{-1}$ admits a unique upper-lower factorisation
\begin{align}\label{C11} &
S^{-1}=(E_T^*)^{-1}   \br E^*\br E      E_T^{-1}
\end{align}
such that
\begin{align}\label{C12} &
 E_T=\sqrt{2}\left(I+\frac{1}{2}\int_0^xs^{\prime}(x-t) \, \cdot \, dt\right), \quad \br E=I+\int_0^x\br \cle(x,t) \, \cdot \, dt, 
\end{align}
where   $\br \cle(x,t)$ is continuous in the triangle $0\leq t\leq x\leq 2\ell$.

Moreover, the operators $E_T^{\pm 1}$  map differentiable functions with a square-integrable derivative on $[0,2 \ell]$
onto differentiable functions with a square-integrable derivative on $[0, 2\ell]$.
\end{Pn}
\begin{Rk}\label{RkV} Since the integral part of $E_T$ $($in \eqref{C12}$)$ is a Volterra operator, the operator $E_T$ has a bounded inverse.
$($On the basic properties of the integral Volterra operators see, for instance, \cite[Chapters 2,9]{Grip}.$)$
In particular, according to Theorems 3.1 and 3.5 in \cite[Chapter 2]{Grip}  the operator $E_T^{-1}$
is again a triangular convolution operator. From the series representation of the Volterra kernel $r(x-t)$ of $E_T^{-1}$
$($see \cite[p. 37]{Grip}$)$ it easily follows that  $r(t)$ is square-integrable.
\end{Rk}
\begin{Rk}\label{RkV2}
Since the kernel of $\br E$ is continuous, it is easy to see that the kernel of $\br E^{-1}$ is continuous in the triangle $0\leq t\leq x\leq 2\ell$ as well
$($see, e.g., \cite[p. 244]{Grip} or consider the series $(I-\Ups)^{-1}=\sum_{i=0}^\infty \Ups^i$ for Volterra operators $\Ups$ of the first kind with continuous kernels$)$.
\end{Rk}
{\it Proof of Proposition \ref{Pn}}.
It is easily proved directly and follows from \cite[Theorem 1.1.3]{SaL3} (see also various references therein)
that $S$ satisfies the operator identity \eqref{C3}, where $A$ is given in \eqref{3.3},  $ \Phi_1, \Phi_2 \in B\big(\BC^{p}, L^2_p(0, \, 2\ell)\big)$
and $\Phi_1 g= \Phi_1(x)g$,  $\Phi_2 g\equiv g$. 

In view of \eqref{C}, the operators $E_T^{\pm 1}$ commute with $A$:
\begin{align}\label{C12+} &
E_T^{\pm 1}A=AE_T^{\pm 1}.
\end{align}
 Hence, using the transformation (operator) $E_T^{-1}$ and the equality $J=\Th j\Th^*$
from \eqref{1.18} we transform \eqref{C3} into the identity
\begin{align}\label{C13} &
\br A\br S-\br S\br A^*=\I \br \Pi j \br\Pi^*,
\end{align}
where $\br \Pi=\begin{bmatrix}\br \Phi_1 & \br \Phi_2\end{bmatrix}$ $\big(\br \Phi_k\in B\big(\BC^p,L^2_p(0,2\ell)\big)\big)$,
\begin{align}\label{C14} &
\br A=-A, \quad \br S=E_T^{-1} S (E_T^{-1})^*, \quad \br \Pi=-E_T^{-1}\Pi j \Th .
\end{align}
The substitution of $\Pi$ with $-\Pi j$ in   the last equality in \eqref{C14} is connected with the substitution of $A$ with $\br A$.
This equality yields
\begin{align}\label{C15} &
\br \Phi_1 g=\br \Phi_1(x)g =\frac{1}{\sqrt{2}}\big(E_T^{-1}(I_p-\Phi_1(x))\big)g, \\
& \label{C16}
 \br \Phi_2 g=\br \Phi_2(x)g =\frac{1}{\sqrt{2}}\big(E_T^{-1}(\Phi_1(x)+I_p)\big)g,
\end{align}
where $E_T^{-1}$ is applied to each column of $I_p\pm\Phi_1(x)$. Similar to the proof of \eqref{3.10}, we obtain
\begin{align}\label{C17} &
\br \Phi_2(x)\equiv I_p.
\end{align}
Indeed, in view of  \eqref{C12}, \eqref{C9} and \eqref{C2+}  we have 
\begin{align}\label{C18} &
E_T I_p=\frac{1}{\sqrt{2}}\big(\Phi_1(x)+I_p\big).
\end{align}
It is easy to see that \eqref{C17} follows from \eqref{C16} and \eqref{C18}.

Similar to the case of $E_0$ in the proof of Lemma \ref{LaSim}, relations \eqref{C12+} and \eqref{C18}
imply that  the operators $E_T^{\pm 1}$  map differentiable functions with a square-integrable derivative on $[0,2 \ell]$
onto differentiable functions with a square-integrable derivative on $[0, 2\ell]$. In particular, the matrix function
$\br \Phi_1(x)$ given by \eqref{C15} is a differentiable function with a square-integrable derivative on $[0, 2\ell]$. 
In view of \eqref{C2+}  and \eqref{C12+}, we also have
\begin{align} &      \label{C19}
\br\Phi_1(0)= 0.
\end{align}

In view of  the properties of $\br \Phi_k(x)$ deduced above and \cite[Remark 1.5]{ALS-HLB},
there is a unique solution $\br S$ of the operator identity (equation) \eqref{C13}, which has the form \cite[(1.9)]{ALS-HLB} considered in \cite[Lemma 3.1]{ALS-HLB}.
From the inequality $S>0$ and the definition of $\br S$ in \eqref{C14}, we obtain $\br S> 0$.
Since $S$ has a bounded inverse, $\br S$ has a bounded inverse as well.
Now, \cite[Lemma 3.1]{ALS-HLB} implies that $\breve S^{-1}$ admits a unique factorisation
\begin{align} &      \label{C19+}
\breve S^{-1}=\br E^*\br E
\end{align}
with $\br E$ as in the statement of our proposition.
The unique factorisation \eqref{C11} follows. Thus, the proposition is proved.

Recall that $S_{\xi}=P_{\xi}S P_{\xi}^*$, where $P_{\xi}$ is defined in \eqref{3.12}, and consider $S_{\xi}^{-1}$. 
\begin{Rk}\label{RkRed} The operators  $\breve E$ and $E_T^{-1}$ are lower triangular and so their product
is lower triangular as well. Now, using  \eqref{C11}, \eqref{C19+} and the considerations from \cite[p. 560]{ALS15}, including the equalities 
which in our case take the form 
\begin{equation}      \label{R}
P_{\xi}\br E=P_{\xi}\br EP_{\xi}^*P_{\xi}, \quad
P_{\xi}E_T^{-1}P_{\xi}^*= (P_{\xi}E_T P_{\xi}^*)^{-1} \quad (\xi<2\ell),
\end{equation}
we similar to \eqref{3.15-} derive
\begin{equation}      \label{C20}
S_{\xi}^{-1}=P_{\xi}(E_T^*)^{-1}   \br E^* P_{\xi}^*P_{\xi}\br E      E_T^{-1}P_{\xi}^*=(E_{T,\,\xi}^*)^{-1}\br E_{\xi}^*\br E_{\xi}E_{T,\,\xi}^{-1}, \quad 
\br S_{\xi}^{-1}=\br E_{\xi}^*\br E_{\xi}, 
\end{equation}
where $\br E_{\xi}=P_{\xi}\br E P_{\xi}^*$, $ \br S_{\xi}=P_{\xi}\br S P_{\xi}^*$, and $\br E_{T,\, \xi}=P_{\xi} E_T P_{\xi}^*$.

In view of \eqref{C10}, \eqref{C9}, $S_{\xi}$ is a convolution operator of the same form as $S$ $($only $\Phi_1(x)$ is considered
on $[0,\xi]$ instead of $[0,2\ell])$.
According to the first equalities in \eqref{C12} and \eqref{C9}, $E_{T,\,\xi}$ has the same form as $E_T$, but $\Phi_1(x)$ is again considered
on $[0,\xi]$ $($instead of $[0,2\ell]$ for $E_T)$. Clearly, the Volterra kernel of $\br E_{\xi}$ is continuous. 

It follows that the first factorisation in \eqref{C20} is the unique factorisation of the type described in Proposition \ref{Pn}.
\end{Rk}

\begin{Rk} \label{RkRed2} Applying $P_{\xi}$ from the left
and $P_{\xi}^*$ from the right to the operator identity \eqref{C13} one can see that $\br S_{\xi}$ satisfies a similar to \eqref{C13}
operator identity $\br A_{\xi}\br S_{\xi}-\br S_{\xi}\br A_{\xi}^*=\I P_{\xi}\br \Pi j \br\Pi^*P_{\xi}^*$, where $\br \Phi_1(x)$ $($and so $\Phi_1(x))$
is considered on $[0,\xi]$ instead of $[0,2\ell]$. Thus, $\br S_{\xi}$ has the same form as $\br S$.
\end{Rk}
Let the operator $S=S_{\tau}$ of the form \eqref{C1} determine matrix functions $\b_{\tau}(x)$ and $\g_{\tau}(x)$ via formulas
\begin{align}\label{C21} &
\b_{\tau}(x)=\br E      E_T^{-1}\begin{bmatrix} \Phi_1(x) & I_p\end{bmatrix},
\\ \label{C22} &
\g_{\tau}(x)= \frac{1}{\sqrt{2}} \left( \begin{bmatrix}-I_p & I_p\end{bmatrix}
- \int_0^x (\Phi_1^{\prime}(\xi))^* S_{x}^{-1}\begin{bmatrix}\Phi_1(\xi) & I_p \end{bmatrix}d\xi \right),
\end{align}
which are similar to the formulas \eqref{3.9} and \eqref{A9} for $\b$ and $\g$, respectively.
\begin{La}\label{LaDif} Let the conditions $(i)-(iii)$ hold. Then, the matrix functions $\b_{\tau}$ and $\g_{\tau}$ given by \eqref{C21} and \eqref{C22}, respectively, 
are differentiable on $[0,2\ell]$, with square-integrable  derivatives. Moreover, we have
\begin{align}\label{C22+} &
\b_{\tau}(0)=\frac{1}{\sqrt{2}}\begin{bmatrix}I_p & I_p\end{bmatrix}, \quad \g_{\tau}(0)=\frac{1}{\sqrt{2}}\begin{bmatrix}-I_p & I_p\end{bmatrix}, 
\\ \label{C23-} &
\b_{\tau}(x)J\g_{\tau}(x)^*\equiv 0, \quad \g_{\tau}^{\prime}(x)J\g_{\tau}(x)^*\equiv 0,
\\ &  \label{C23!}
\b_{\tau}^{\prime}(x)J\b_{\tau}(x)^*\equiv 0.
\end{align}
\end{La}
\begin{proof} 
The first equality in \eqref{C22+} follows from \eqref{C2+}, \eqref{C12}, and \eqref{C21}. The second equality in \eqref{C22+} is immediate from \eqref{C22}.

According to Proposition \ref{Pn}, the matrix function $E_T^{-1}\begin{bmatrix} \Phi_1(x) & I_p\end{bmatrix}$ is differentiable on $[0,2\ell]$,
with square-integrable  derivative. Moreover, the Volterra kernel of $\br E$ is continuous, which implies that $\b_{\tau}(x)$ is continuous (at least).
Taking into account \eqref{C20}, we factorise $S_x^{-1}$ in \eqref{C22} and differentiate both parts of
\eqref{C22}:
\begin{align} \label{C23} &
\g_{\tau}^{\prime}(x)= -\frac{1}{\sqrt{2}} \big(\br E      E_T^{-1}\Phi_1^{\prime}\big)(x)^*
\big(\br E      E_T^{-1}\begin{bmatrix}\Phi_1& I_p \end{bmatrix}\big)(x).
\end{align}
Using \eqref{C21} we rewrite \eqref{C23} as
\begin{align} \label{C24} &
\g_{\tau}^{\prime}(x)= -\frac{1}{\sqrt{2}} \big(\br E      E_T^{-1}\Phi_1^{\prime}\big)(x)^*\b_{\tau}(x).
\end{align}
Thus,  $\g_{\tau}$ is differentiable on $[0,2\ell]$, with square-integrable  derivative.

In order to prove the first equality in \eqref{C23-}, we recall that the conditions $(i)-(iii)$ yield the operator identity
\eqref{C3}. Note also that, according to \eqref{C11}, $S$ in \eqref{C3} admits factorisation $S=E_T \br E^{-1}(E_T \br E^{-1})^*$.
Thus, multiplying \eqref{C3} by $\br E E_T^{-1}$ from the left and by $(\br E E_T^{-1})^*$ from the right and taking into account
that $\Pi=\begin{bmatrix} \Phi_1 & \Phi_2\end{bmatrix}$, $\Phi_1g=\Phi_1(x)g, \,\, \Phi_2 g\equiv g$, and the definition 
\eqref{C21} of $\b_{\tau}$ is valid, we obtain
\begin{align}\label{C25} &
\br E E_T^{-1}A E_T \br E^{-1}-(\br E E_T^{-1}A E_T \br E^{-1})^*=\I \b_{\tau}(x)J\int_0^{2\ell}\b_{\tau}(t)^*\, \cdot \, dt.
\end{align}
Since $\br E E_T^{-1}A E_T \br E^{-1}$ is a lower  triangular integral operator and its adjoint is an upper triangular integral operator,
formula \eqref{C25} implies that
\begin{align}\label{C26} &
E_{\tau}A E_{\tau}^{-1}=\I \b_{\tau}(x)J\int_0^{x}\b_{\tau}(t)^*\, \cdot \, dt, \quad E_{\tau}:=\br E E_T^{-1}.
\end{align}
Similar to $\b$, we partition $\b_{\tau}$ into two $p\times p$ blocks $\b_{\tau}(x)=\begin{bmatrix}(\b_{\tau})_1(x) & (\b_{\tau})_2(x) \end{bmatrix}$.
Recall that $A=\I \int_0^x\cdot \, dt$. Therefore, it follows from \eqref{C2+} and \eqref{C21} that
\begin{align}\nn
\big(E_{\tau}A E_{\tau}^{-1}\big)\big(E_{\tau} \Phi_1^{\prime}\big)(x)& =\big(E_{\tau} A\Phi_1^{\prime}\big)(x)=\I E_{\tau} \big(\Phi_1(x)-I_p\big)
\\ \label{C27} &
=\I\big((\b_{\tau})_1(x) - (\b_{\tau})_2(x)\big),
\end{align}
where operators are as usual applied to the columns of the matrix functions. (Compare \eqref{C27} with \eqref{A10}.) Using  \eqref{C26}, we rewrite \eqref{C27} as
\begin{align}
\label{C28} &
(\b_{\tau})_1(x) - (\b_{\tau})_2(x)=\b_{\tau}(x)J\int_0^x\b_{\tau}(t)^*\big(E_{\tau} \Phi_1^{\prime}\big)(t)dt.
\end{align}
In view of  the equality $S_{x}^{-1}=(E_{T,\,x}^*)^{-1}\br E_{x}^*\br E_{x}E_{T,\,x}^{-1}$ and relations \eqref{R}, \eqref{C21} and \eqref{C22}, equality \eqref{C28} takes the form
\begin{align}\nn
(\b_{\tau})_1(x) - (\b_{\tau})_2(x)&=\b_{\tau}(x)J\int_0^x\big(\br E_{x} E_{T,x}^{-1}\begin{bmatrix}\Phi_1(t) & I_p \end{bmatrix}\big)^*\big(\br E_{x} E_{T,x}^{-1} \Phi_1^{\prime}\big)(t)dt
\\ & \label{C29} 
=\b_{\tau}(x)J\int_0^x\big(S_x^{-1}\begin{bmatrix}\Phi_1(t) & I_p \end{bmatrix}\big)^* \Phi_1^{\prime}(t)dt
\\ & \nn
=\b_{\tau}(x)J\big(\begin{bmatrix}-I_p & I_p\end{bmatrix}-\sqrt{2}\g_{\tau}(x)\big)^*.
\end{align}
The  first equality in \eqref{C23-}, that is, the equality $\b_{\tau}(x)J\g_{\tau}(x)^*\equiv 0$  immediately follows from \eqref{C29}. Furthermore, the  first equality in \eqref{C23-} and formula
\eqref{C24} imply the second equality in \eqref{C23-}.

It remains to prove that $\b_{\tau}$ is differentiable and that \eqref{C23!} holds. First, we show that 
\begin{align}\label{C30} &
\b_{\tau}(x)J\b_{\tau}(x)^*\equiv I_p.
\end{align}
For this purpose, we note that  $E_{\tau}^{\pm 1}$ in \eqref{C26} are integral Volterra operators and may be represented in the form
\begin{align}\label{C31} &
E_{\tau}^{-1}=\sqrt{2}\left(I+\int_0^x\G_{\tau}(x,t) \, \cdot \, dt\right), \quad E_{\tau}=\frac{1}{\sqrt{2}}\left(I+\int_0^x\cle_{\tau}(x,t) \, \cdot \, dt\right),
\end{align}
where $\G_{\tau}$ and $\cle_{\tau}$ are Hilbert-Schmidt kernels (for further details on these kernels see Remarks \ref{RkV} and \ref{RkV2}). 
Next, we rewrite the similarity relation in \eqref{C26} in terms of the equality of the corresponding
integral operator kernels:
\begin{align}\nn
\b_{\tau}(x)J\b_{\tau}(\xi)^*=& I_p+\int_{\xi}^x\big(\cle_{\tau}(x,t) +\G_{\tau}(t,\xi) \big)dt
\\ \label{C32}  &
+\int_{\xi}^x\int_\xi^{t_2}\cle_{\tau}(x,t_2)\G_{\tau}(t_1,\xi)dt_1dt_2.
\end{align}
We derive \eqref{C30} setting $\xi=x$ in \eqref{C32}.

In view of  the second relations in \eqref{C22+} and \eqref{C23-}, the matrix function $\g=\g_{\tau}$ satisfies the conditions
of Proposition \ref{PnBG}. Let us consider the matrix function $\b(x)$ recovered in Proposition \ref{PnBG}  from $\g(x)$ and show
that $\b=\b_{\tau}$. Note that $\b$ satisfies \eqref{A20-} and $ \b J \b^*\equiv I_p$ follows. 
Since $ \b J \g_{\tau}^*\equiv \b_{\tau}  J \g_{\tau}^*\equiv 0$ and  $ \b J \b^*\equiv \b_{\tau}  J \b_{\tau}^*\equiv I_p$,
we derive
\begin{align}\label{C33} &
\b_{\tau}(x)=\om(x)\b(x), \quad \om(x)^*=\om(x)^{-1},
\end{align}
and we will show that $\om(x) \equiv I_p$.

Indeed, formula \eqref{C26} may be rewritten in the form
\begin{align}\label{C26'} &
E_{\tau}A =\I \b_{\tau}(x)J\int_0^{x}\b_{\tau}(t)^*\, \cdot \, dt \, E_{\tau}.
\end{align}
Similar to  Lemma \ref{LaSim}, which follows from Corollary \ref{CySim} without specifically using Dirac systems, we have
\begin{equation}
\label{C34} \I \b(x)J \int_0^x \b(t)^*\, \cdot \,dt= E A  E^{-1},  \end{equation}
where $E$ has the same properties as $E$ in Lemma \ref{LaSim}.
Now, relations \eqref{C33} and \eqref{C34} yield
\begin{align}\label{C35} &
E_{\om}A= \I \b_{\tau}(x)J \int_0^x \b_{\tau}(t)^*\, \cdot \,dt \, E_{\om}, \quad E_{\om}f:=\om(x)Ef \quad (f\in L^2_p(0,2\ell)).
\end{align}
According to \eqref{3.4}, \eqref{C33} and to \eqref{C21} we have
\begin{align}\label{C36} &
E_{\om}I_p=E_{\tau} I_p=(\b_{\tau})_2(x).
\end{align}
Relations \eqref{C26'}, \eqref{C35} and \eqref{C36} yield
\begin{align}\label{C36+} &
E_{\om} A^k I_p = E_{\tau}A^k I_p  \quad (k\geq 0).
\end{align}
Finally, it is easy to see that
\begin{align}\label{C37} &
\overline{\spa}\bigcup_{i=0}^{\infty}\{A^ig: \, g\in \BC^p\}=L^2_p(0,2\ell),
\end{align}
where $\overline{\spa}$ is the closure of the linear span. Taking into account \eqref{C36+} and \eqref{C37}, we obtain
\begin{align}\label{C38} &
E_{\om}=E_{\tau}.
\end{align}
Comparing representations of $E$ and $E_{\tau}$ in \eqref{3.3} and \eqref{C31}, respectively, and recalling
that $ E_{\om}=\om(x)E$, we deduce that $\om(x)f (x)=f(x)$ for the functions $f\in L^2_p(0,2\ell)$, that is, $\om(x)\equiv I_p$. Thus,
\begin{align}\label{C39} &
\b_{\tau}(x)=\b(x),
\end{align}
which implies that $\b_{\tau}$ is
differentiable on $[0,2\ell]$, with square-integrable  derivatives. In view of \eqref{C39} and the second relation in \eqref{A20-}, the 
identity \eqref{C23!} holds.
\end{proof}
\noindent {\bf 4.} {\it Proof of Theorem \ref{TmSets}}.  

The necessity of the properties $(i)-(iii)$ was already proved earlier (see paragraph {\bf 1} of this section).
Here, we show their sufficiency.

Step 1. Let us consider Dirac system on $[0,\ell]$ with the potential
\begin{align}\label{C40} &
v(x/2)=2\I \b_{\tau}^{\prime}(x)J\g_{\tau}(x)^*,
\end{align}
where $\b_{\tau}$ and $\g_{\tau}$ are constructed above. We will show that $\tau$ $\big($under conditions $(i)-(iii)\big)$ is
a spectral function of this system. Since $\b_{\tau}J\g_{\tau}^*\equiv 0$, we also have 
$$\g_{\tau}^{\prime}(x)J\b_{\tau}(x)^*+\big(\b_{\tau}^{\prime}(x)J\g_{\tau}(x)^*\big)^*\equiv 0,$$
which together with \eqref{C40} yields
\begin{align}\label{C41} &
\g_{\tau}^{\prime}(x)J\b_{\tau}(x)^*=-\big(\b_{\tau}^{\prime}(x)J\g_{\tau}(x)^*\big)^*=-\frac{\I}{2}v(x/2)^*.
\end{align}

We introduce the matrix function
\begin{align}\label{C42} &
y(x)=\begin{bmatrix}\b_{\tau}(2x)\\ \g_{\tau}(2 x) \end{bmatrix}J\Th j,
\end{align}
and recall relations \eqref{C22+}--\eqref{C23!} which imply
\begin{align}\label{C43} &
\b_{\tau}(x)J\g_{\tau}(x)^*\equiv 0, \quad \b_{\tau}(x)J\b_{\tau}(x)^*\equiv I_p, \quad \g_{\tau}(x)J\g_{\tau}(x)^*\equiv -I_p.
\end{align}
Taking into account the last equality in \eqref{1.18} (i.e., the equality $\Th j \Th^*=J$) and formulas \eqref{C42}, \eqref{C43}, we obtain
\begin{align}\label{C44} &
y(x)j y(x)^*\equiv j.
\end{align}
From \eqref{C42} and \eqref{C44} we derive
\begin{align}\label{C45} 
y^{\prime}(x)&=2\begin{bmatrix}\b_{\tau}^{\prime}(2x)\\ \g_{\tau}^{\prime}(2x) \end{bmatrix}J\Th j j y(x)^* j y(x).
\end{align}
Hence, it follows from \eqref{C23-}, \eqref{C23!} and  \eqref{C40}--\eqref{C42} that
\begin{align}
 &\label{C46} 
y^{\prime}(x)=2\begin{bmatrix}\b_{\tau}^{\prime}(2x)\\ \g_{\tau}^{\prime}(2x) \end{bmatrix}J \begin{bmatrix}\b_{\tau}(2x)^* & \g_{\tau}(2x)^* \end{bmatrix}j y(x)=\I j V(x) y(x),
\end{align}
where $V$ has the form \eqref{1.2}. Moreover, in view of \eqref{1.18}, \eqref{C22+}, and \eqref{C42} we have $y(0)=I_{2p}$.
Thus,
\begin{align}\label{C47} &
y(x)=u(x,0),
\end{align}
where $u(x,z)$ is the normalized fundamental solution of the Dirac system with the potential given by \eqref{C40}.
According to \eqref{C42} and \eqref{C47}, the equality $u(x/2,0)\Th^*=\begin{bmatrix}\b_{\tau}(x) \\ \g_{\tau}(x)\end{bmatrix}$ is valid. Therefore,
$\b$ and $\g$ corresponding to Dirac system with $v$ of the form \eqref{C40} via \eqref{1.15} and \eqref{1.16} satisfy relations
\begin{align}\label{C48} &
\b(x)=\b_{\tau}(x), \quad \g(x)=\g_{\tau}(x).
\end{align}

Step 2. The operator $E$ in Section \ref{inv} satisfies relations \eqref{3.2}--\eqref{3.4}. In view of  \eqref{C21}, \eqref{C26}, \eqref{C31}, and \eqref{C48}
$E_{\tau}$ satisfies the same relations. Hence, similar to the equality  \eqref{C38}, we obtain for this $E$ the equality
\begin{align}\label{C49} &
E=E_{\tau}.
\end{align}
In view of \eqref{3.9}, \eqref{C21}, \eqref{C48}, and \eqref{C49}, $\Phi_1(x)$ generated by $\tau$ (which we further  denote by $(\Phi_{\tau})_1(x)$
coincides with $\Phi_1(x)$  corresponding to Dirac system with $v$ of the form \eqref{C40}:
\begin{align}\label{C50} &
\Phi_1(x)=(\Phi_{\tau})_1(x).
\end{align}
Relations \eqref{3.8}, \eqref{C11} and \eqref{C49} imply $S=S_{\tau}$. Hence, taking into acount  \eqref{C50}, we see that the transfer
matrix function $w_A(2\ell, z)$ corresponding to the constructed above Dirac system coincides with the transfer matrix function 
$$(w_A)_{\tau}(z)=I_{2p}+\I z J\Pi_{\tau}^*S_{\tau}^{-1}(I-zA)^{-1}\Pi_{\tau} \quad \big(\Pi_{\tau}:=\begin{bmatrix}(\Phi_{\tau})_1 &\Phi_2\end{bmatrix}\big),$$
generated by $\tau$. Remark \ref{RkSpF} and relations \eqref{3.17}, \eqref{3.18}, and \eqref{3.20} show that the matrix functions
(M\"obius transformations) from the class $\cln(\A)$ determined by $w_A(2\ell, z)$ describe Weyl functions of the constructed Dirac system.
Moreover, in view of Remark \ref{RkSpF}, distribution functions $\tau$ in Herglotz representations \eqref{1.13} of the functions from $\cln(\A)$ give us the set of spectral functions.

According to \cite[Theorem 1.27]{SaSaR} (see also \cite[Chapter 4]{SaL2}) and to the equality $w_A(2\ell, z)=(w_A)_{\tau}(z)$, the set of functions $\tau$ satisfying condition $(i)$ and representations
\eqref{C1} and \eqref{C2}, \eqref{C2-} coincides with the set of $\tau$ obtained using $\cln(\A)$ above, and so with the set of spectral functions. 
Thus, the matrix function $\tau$ with properties $(i)-(iii)$, which we consider in the theorem, is one of the
spectral functions of the constructed Dirac system with $v$ of the form \eqref{C40}. The theorem is proved.

Recalling Definition \ref{DnSF}, we obtain the following corollary of Theorem~\ref{TmSets}.
\begin{Cy} \label{Cysa} Let a given  $p\times p$ matrix function $\tau(t)$ on $\BR$ have the properties $(i)$. 
Assume that $\tau$ satisfies conditions $(ii)$ and $(iii)$ for some values $\ell=\ell_k>0$ $(1 \leq k<\infty)$,
where the sequence $\{\ell_k\}$ tends to infinity for $k$ tending to infinity.
Then,  $\tau(t)$ is the spectral function of  some Dirac system on $[0,\infty)$ with a locally square-integrable 
potential.
\end{Cy}
 \section{$PW$-sampling measures and reformulations} \label{PW}
\setcounter{equation}{0} 
\subsection{Reformulation of the conditions $(ii)$ and $(iii)$}
Let us recall some arguments from paragraph {\bf 1} of Section \ref{set}. According to Remark \ref{RkBST}, the weak convergence of the integral on the right-hand side of \eqref{C1} 
suffices for the boundedness of $S$. If we have this   weak convergence and if the condition $(i)$ 
is satisfied, the integral in \eqref{C2-} (i.e., the operator $\Lam$) is well defined. 
Therefore, the continuous matrix function $\Phi_1(x)$ is also well defined.

Now, under condition $\nu=\nu^*$ (see \eqref{C2++}),  it follows from the proof of Lemma~\ref{LaId} that the operator identity \eqref{C3} holds.
Moreover, we have \eqref{C2+}, that is, $\Phi_1(0)=I_p$. Hence, the bounded solution $S$ of the equation (operator  identity) \eqref{C3} is given by
\eqref{C10}. Here, we used \cite[Corollary D.4]{SaSaR} (see also \cite[Theorem 1.1.3]{SaL3}) and the differentiability of $\Phi_1(x)$.

It follows from the representation \eqref{C1} that $S=S^*\geq 0$. Since $\Phi_1^{\prime}(x)$ is square-integrable, the integral part
of the operator $S$ (on the right-hand side of \eqref{C10}) is a Hilbert-Schmidt operator. Thus, the existence and boundedness
of $S^{-1}$ is equivalent to the condition $Sf=0$ iff $f=0$, that is, ${\mathrm{Ker}}\, S=0$. We derived the following proposition.
\begin{Pn}\label{PnReform} Under condition $(i)$ $($see paragraph {\bf 1} of Section \ref{set}$)$, conditions $($properties$)$
$(ii)$ and $(iii)$ are equivalent to the conditions\\
$(II)$ The integral on the right-hand side of \eqref{C1} weakly converges.\\
$(III)$ The matrix function $\Phi_1(x)$ determined by the operator $\Phi_1$ of the form \eqref{C2} is differentiable
and the derivative $\Phi_1^{\prime}(x)$ is square-integrable on $(0,2\ell)$.\\
$(IV)$ The matrix $\nu$ of the form \eqref{3.24} is self-adjoint and we have
$${\mathrm{Ker}}\, S=0 \quad {\mathrm{for}} \quad S=2I+\int_0^{2\ell}s^{\prime}(x-t) \, \cdot \, dt,$$
where $s(x)=\Phi_{1}(x)$ for $x>0$ and $s(x)=-s(-x)^*$ for $x<0$. 
\end{Pn}

\subsection{PW-sampling measures}
Assuming that the right-hand side of \eqref{C1} weakly converges and taking into account \eqref{3.23}, we see
that the operator $S=S_{2\ell}$ in \eqref{C1} satisfies the relations
\begin{align}\nn
\langle S f,f\rangle_{L^2}&=\int_{-\infty}^{\infty}\left(\int_0^{2\ell}\E^{-\I \eta_1 t}f(\eta_1)d\eta_1\right)^*d\tau(t)\int_0^{2\ell}\E^{-\I \eta_2 t}f(\eta_2)d\eta_2
\\  \label{P1} &
\\ \nn &
=\int_{-\infty}^{\infty}\left(\int_{-\ell}^{\ell}\E^{-\I \xi_1 t}f(\xi_1+\ell)d\xi_1\right)^*d\tau(t)\int_{-\ell}^{\ell}\E^{-\I \xi_2 t}f(\xi_2+\ell)d\xi_2,
\end{align}
where $f\in L^2_p(0,2\ell)$ and $\xi_k=\eta_k-\ell$ $(k=1,2)$.

Let us consider in greater detail the case 
\begin{align}\label{P2} &
p=1,
\end{align}
that is, the case of scalar functions $\tau$.  This case is related to the case of  canonical systems \eqref{1.14} with $2\times 2$
Hamiltonians $H(x)$,
which is studied much better than the no less important  (but more difficult) case of  $H(x)$ of the order $2p$ ($p>1$). Fundamental results on 
the case $p=1$ have been obtained by L.~de Branges and by M.G.  Krein, see also \cite{LW, MP, Su} for further references.

Two results on $PW$-sampling measures (see \cite[Subsection~3.1]{MP}) are essential for us in the case $p=1$.
(Note that we use the definitions and formulations from \cite{MP} which may somewhat differ from other sources.)
We will need the following definitions.
\begin{Dn}\label{DnPW} By $PW_{\ell}$ is denoted a subclass of functions  $h(t)$ from $\break L^2(-\infty,\infty)$ which admit a representation
\begin{align}\label{P2+} &
h(t)=\int_{-\ell}^{\ell}\E^{-\I t\xi}f(\xi)d\xi, \quad f \in L^2(-\ell,\ell).
\end{align}
\end{Dn}
Clearly, the functions $h(t)$ have entire extensions (and PW in the notation stands for the Paley-Wiener space).
\begin{Dn}\label{DnPWS} A  positive measure $d\tau$ on $\BR$
is $PW_{\ell}$-sampling if there is some $C>0$ such that
\begin{align}\label{P3} &
C^{-1}\|h\|_{L^2}\leq \|h\|_{L^2(d\tau)}\leq C \|h\|_{L^2} \quad {\mathrm{for \,\, any}} \quad h\in PW_{\ell},
\end{align}
where $\|\cdot\|_{L^2}$ is the norm in $L^2(-\infty,\infty)$ and  $\|\cdot \|_{L^2(d\tau)}$ is the norm in $L^2(d\tau)$.

A measure $d\tau$ is $PW$-sampling if $d\tau$ is $PW_{\ell}$-sampling for all $\ell>0$.
\end{Dn}
\begin{Dn}\label{DnSeq} A real sequence $\{\la_k\}_{k\in \BZ}$ is sampling for $PW_{\ell}$ if it is separated, that is,
$|\la_i-\la_k|>\ve>0$ $($for some $\ve$ and all $i\not=k)$ and there is some $C>0$ such that
$$C^{-1}\|h\|_{L^2}^2\leq {\sum_{k=-\infty}^{\infty}|h(\la_k)|^2}\leq C \|h\|_{L^2}^2 \quad {\mathrm{for \,\, any}} \quad h\in PW_{\ell}.$$
\end{Dn}
Now, we formulate an important propositions on $\tau$.
\begin{Pn} \label{PnOrS}  \cite{OrS} A positive measure $d\tau$ on $\BR$ is $PW_{\ell}$-sampling   iff
$\tau(x+1-0)-\tau(x+0)$ is uniformly bounded on $\BR$ and  there exist a sampling $($for $PW_{\ell})$ sequence $\{\la_k\}_{k\in \BZ}$
and a number $\de>0$ such that $\displaystyle{\int_{U_k}d\tau>\de}$ $\break (-\infty<k<\infty)$ for some disjoint neighborhoods $U_k$ of $\la_k$.
\end{Pn} 
For the case of the $PW$-sampling measures $d\tau$ we need one definition more.
\begin{Dn}\label{DnCap} An interval $\cli\in \BR$ is called $\de$-massive $(\de>0)$ with respect to a positive measure $d\tau$ if
$|\cli|\geq \de$ and $\int_{\cli}d\tau \geq \de$. 

The $\de$-capacity $C_{\de}(\cli)$ $($for a positive measure $d\tau$ and an interval $\cli)$ is the maximal number of disjoint $\de$-massive intervals intersecting $\cli$.
\end{Dn}
\begin{Pn} \label{PnMP}  \cite{MP} A positive measure $d\tau$ on $\BR$ is $PW$-sampling   iff
$\tau(x+1-0)-\tau(x+0)$ is uniformly bounded on $\BR$ and  for any $r>0$ there exist $c>0$ and $\de>0$
such that $C_{\de}(\cli)\geq r|\cli|$ for all $\cli$ satisfying the inequality $|\cli|\geq c$.
\end{Pn}

In view of \eqref{P1} and Proposition \ref{PnOrS}, the important condition  $(ii)$ in the paragraph {\bf 1} of Section \ref{set}
may be reformulated, that is, we have the following corollary.
\begin{Cy}\label{Cyii} Let $p=1$. Then, the conditions (ii) on the boundedness of the operator $S$ $($given by
\eqref{C1}$)$ and of its inverse may be substituted by the equivalent conditions:\\
$(ii^{\prime})$ For the given  positive measure $d\tau$, the difference $\tau(x+1-0)-\tau(x+0)$ is uniformly bounded on $\BR$.  Moreover, for this $d\tau$,
there exist a sampling $($for $PW_{\ell})$ sequence $\{\la_k\}_{k\in \BZ}$
and a number $\de>0$ such that $\displaystyle{\int_{U_k}d\tau>\de}$ $\break (-\infty<k<\infty)$ for some disjoint neighborhoods $U_k$ of $\la_k$.
\end{Cy}
We also have a corollary of Proposition \ref{PnMP} on the condition $(ii)$ in Corollary \ref{Cysa} (for a sequence $\ell_k\to \infty$ or, equivalently,
for all $\ell>0$).
\begin{Cy}\label{Cyiisa} Let $p=1$. Then, the conditions of the boundedness of the operators $S$ $($given by
\eqref{C1}$)$  and of their inverses for all $\ell>0$ may be substituted by the equivalent conditions:\\
$(ii^{\prime})$ For the given  positive measure $d\tau$, the difference $\tau(x+1-0)-\tau(x+0)$ is uniformly bounded on $\BR$.  Moreover, for this $d\tau$ and  
for any $r>0$, there exist $c>0$ and $\de>0$
such that $C_{\de}(\cli)\geq r|\cli|$ for all $\cli$ satisfying the inequality $|\cli|\geq c$.
\end{Cy}
It would be very interesting to obtain analogs of Propositions \ref{PnOrS} and \ref{PnMP} for the matrix valued measures (distribution functions).
In order to formulate these open problems, let functions $f$ in \eqref{P2+} belong to the class $L^2_p(-\ell,\ell)$ and let $\|h\|_{L^2}$ 
(for the functions $h$ corresponding to $f$ in \eqref{P2+}) be the norm of $h$
in $L^2_p(-\infty,\infty)$. 

{\bf Open Problems}\\ 
a) For the  nondecreasing on $\BR$ distribution $p\times p$ matrix functions $\tau$, find conditions under which
the inequalities in \eqref{P3} hold $($for some $C>0)$.\\
b) For the  nondecreasing on $\BR$ distribution $p\times p$ matrix functions $\tau$, find conditions under which
the inequalities in \eqref{P3} hold for all $\ell>0$ (and for some $C=C(\ell)>0$ at each $\ell$). 
\begin{Rk}\label{RkMikh} It is easy to see $($and was communicated once to the author by I.V. Mikhailova$)$ that
the boundedness of the operator $S_{2\ell}$ for some $\ell>0$ $($i.e., the right inequality in \eqref{P3}$)$ yields the
boundedness of the operators $S_{2\ell}$ for all $\ell>0$.
\end{Rk}
Indeed, let  the  operator $S_{2\ell}$ (for some  $\ell>0$)  be bounded. Then, \eqref{P1} yields that the norm of $S_{2r}$ $(r<\ell)$ is bounded
(by the same constant). Therefore, it suffices to show that the norms of $S_{2k\ell}$ $(k\in \BN)$ are bounded. In the case of $S_{2k\ell}$, the integrals
from $0$ till $2\ell$ in \eqref{P1} are substituted by the integrals
from $0$ till $2k\ell$.
We split the interval $[0,2k\ell]$ into $k$ intervals of the length $2\ell$ and the integrals from $0$
till $2k\ell$  into the sums of $k$ corresponding integrals.
In this way,  formula \eqref{P1}  and the inequality $\langle S_{2\ell} f,f\rangle_{L^2}\leq M \langle  f,f\rangle_{L^2}$ 
imply that $\langle S_{2k\ell} f,f\rangle_{L^2}\leq k^2 M \langle  f,f\rangle_{L^2}$, where  $f\in L^2_p(0,2k\ell)$.
We note that the norm of $S_{2\ell}$ (or its upper bound $2\pi C^2$ in the notations from \eqref{P3}) may be taken as $M$
in the inequalities above. 


\appendix
\section{Recovery of $\b$ and $\g$} \label{app}
\setcounter{equation}{0}
Operators $S$ and $\Pi$ for Dirac systems on $[0,\ell]$ with square-integrable potentials were constructed in the paragraph {\bf 2} of Section  \ref{inv}.
Given operators $S$ and $\Pi=\begin{bmatrix}\Phi_1 & \Phi_2\end{bmatrix}$ as in  formulas \eqref{3.8}--\eqref{3.10}, we need to recover
matrix functions $\b(x)$ and $\g(x)$ (introduced in \eqref{1.15} and \eqref{1.16}) in order to solve the inverse problem. This can be done in several
ways (see \cite{SaA3} and related cases in \cite{ALS15} and in \cite[pp. 14--17]{SaSaR}).

{\bf 1.} For instance, using \eqref{3.15} we obtain matrix function
\begin{align}\label{A1} &
H(x)=\b(x)^*\b(x).
\end{align}
It follows from $\b J\b^*\equiv I_p$ that $\det \big(\b_k(x)\big)\not=0$ $(k=1,2)$, and we set 
\begin{align}\label{A2} &
\th:=\big(\b_2(x)^*\big)^{-1}, \quad \wt \b(x):=\begin{bmatrix}0 & I_p \end{bmatrix}H(x).
\end{align}
According to \eqref{A1} and \eqref{A2}, we have
\begin{align}\label{A3} &
\b(x)=\th(x)\wt\b(x).
\end{align}
\begin{Pn}\label{PnBeta} The matrix function $\b(x)$ is recovered using \eqref{A3}, where $\wt \b(x)$ is given by  \eqref{A2},
and the $p \times p$ matrix function $\th(x)$ is the unique solution of the first order differential equation
\begin{align}\label{A4} &
\th^{\prime}(x)+\th(x)\big(\wt\b^{\prime}(x)J\wt \b(x)^*\wt \b_2(x)^{-1}\big)=0, \quad \th(0)=\sqrt{2}I_p.
\end{align}
\end{Pn}
\begin{proof}  Taking into account \eqref{1.21} and \eqref{A3}, we rewrite the equality $\break \b^{\prime}J\b^* \equiv 0$ in \eqref{1.23+} in the form
\begin{align}\label{A5} &
\b^{\prime}J\b^*=\th^{\prime}\th^{-1}\b J\b^*+\th \wt \b^{\,\prime}J\wt \b^{\,*}\th^*=\th^{\prime}\th^{-1}+\th \wt \b^{\,\prime}J\wt \b{\,^*}\th^*=0.
\end{align}
According to \eqref{A1} and \eqref{A2}, we have $\wt \b_2^{-1}=\th^*\th$.  Hence, the equation \eqref{A4} follows from \eqref{A5}. (Note that the initial
condition in \eqref{A4} is immediate from the value of $\b(0)$ in \eqref{1.21}.)
\end{proof}
Next, the matrix function $\g$ is recovered from $\b$. In view of the equality $\b J \g^*\equiv 0$, we represent $\g$ in the form
\begin{align}\label{A6} &
\g(x)=\vt(x)\wt \g(x), \quad \wt \g(x):=\frac{1}{2}\begin{bmatrix} -\big(\b_2(x)^*\big)^{-1}& \big(\b_1(x)^*\big)^{-1} \end{bmatrix},
\end{align}
where $\b(x)J\wt \g(x)^*\equiv 0$ and $\vt(x)$ is a $p\times p$ matrix function.
\begin{Pn}\label{PnGamma} The matrix function $\g(x)$ is recovered from $\b(x)$ using  \eqref{A6}, where
$\vt(x)$ is the unique solution of the first order differential equation
\begin{equation}\label{A7}
\vartheta^{\prime} (x)=- \vartheta (x)  \big(\tilde \g^{\prime}(x)
  J \tilde \g(x)^{*}\big)\big(\tilde \g(x) J \tilde \g(x)^*\big)^{-1},
\qquad \vartheta(0)=I_p.
\end{equation}
\end{Pn}
\begin{proof}  The proof is similar to the proof of  Proposition \ref{PnBeta}. This time, we rewrite the equality $ \g^{\prime}J\g^* \equiv 0$ in \eqref{1.23+}
using the equality $\g J\g^*\equiv -I_p$ from \eqref{1.22}:
$$ \g^{\prime}  J \g^{*}=
 \vartheta^{\prime}
\vartheta ^{-1} \g J \g^{*}+
\vartheta  \tilde \g^{\prime}  J \g^{*}
=
- \vartheta^{\prime}
\vartheta^{-1}
+ \vartheta  \tilde \g^{\prime}  J \tilde \g^{*} \vartheta^*=0.$$
It also follows from $\g J\g^*\equiv -I_p$ that $\vt^*\vt=-\big(\wt \g J \wt \g^*\big)^{-1}$. Therefore, the formula above
implies \eqref{A7}.
\end{proof}
\begin{Rk} \label{RkG}
The matrix function $\gamma(x)$ is uniquely recovered from $\b$ $($in the proof of Proposition \ref{PnGamma}$)$ using its three properties only:
\begin{align}\label{A8} &
\g(0)=\frac{1}{\sqrt{2}}\begin{bmatrix}-I_p & I_p \end{bmatrix}, \quad \b(x)J \g(x)^*\equiv 0, \quad \g^{\prime}(x)J\g(x)^* \equiv 0.
\end{align}
The property $\g J\g^*\equiv -I_p$  follows from $ \g^{\prime}  J \g^{*}\equiv 0$ and from the initial value $\g(0)$.
\end{Rk}
{\bf 2.}  There is a direct way to recover $\g$ from $S$ and $\Pi$ as well.
\begin{Pn}\label{PnDg} Given a Dirac system on $[0,\ell]$ and the corresponding operators $S$ and $\Pi$, the following
equality holds for $\g$:
\begin{equation} \label{A9}
\g(x)= \frac{1}{\sqrt{2}} \left( [-I_p \quad I_p]
- \int_0^x (\Phi_1^{\prime}(\xi))^* S_{x}^{-1}\begin{bmatrix}\Phi_1(\xi) & I_p \end{bmatrix}d\xi \right),
\end{equation}
where $S_{x}^{-1}$ is applied to the matrix function $\begin{bmatrix}\Phi_1(\xi) & I_p \end{bmatrix}$ columnwise.
\end{Pn}
\begin{proof} Denote the right-hand side of  \eqref{A9} by $\widehat \g$. According to Remark \ref{RkG}, it suffices
to show that relations \eqref{A8} are valid after substitution of $\widehat \g$ (instead of $\g$) in order to prove \eqref{A9}. 

It is immediate that the first equality in \eqref{A8}  holds for $\wh g(0)$.
From \eqref{3.9}--\eqref{3.10+} we derive
\begin{align}  \nn
 \b_1(x) - \b_2(x)&=
E\big(\Phi_1(x)-I_p\big)=E\big(\Phi_1(x)-\Phi_1(0)\big)
\\ \label{A10} &
=-\I \big(EAE^{-1}\big)E \Phi_1^{\prime}(x).
\end{align}
Using \eqref{3.2}, \eqref{3.15-}, \eqref{3.14+} and the equality $\big(P_{x}\b(\xi)\big)^*=\big(E_{x} \begin{bmatrix} \Phi_1(\xi) & I_p \end{bmatrix} \big)^*$
for $\xi\leq x$, we rewrite \eqref{A10} in the form
\begin{align} \nn
\b(x)
\left[ \begin{array}{c} I_p \\ -I_p \end{array} \right] &=
\b(x) J \int_0^x \left[ \begin{array}{c} \Phi_1(\xi)^* \\ I_p \end{array} \right]
E_x^*E_x\Phi_1^{\prime}(\xi)d\xi
\\ &  \label{A11}
=\b(x) J \int_0^x \left[ \begin{array}{c} \Phi_1(\xi)^* \\ I_p \end{array} \right]
S_x^{-1}\Phi_1^{\prime}(\xi)d\xi.
\end{align}
The equality \eqref{A11} may be rewritten as $\b(x)J\widehat \g(x)^*=0$. Thus, $\widehat \g$ satisfies the second relation in \eqref{A8}.

Since $E_{x} \begin{bmatrix} \Phi_1(\xi) & I_p \end{bmatrix} =\b(\xi)$ for $\xi\leq x$,
taking into account the definition of $\widehat \g$ and the formula \eqref{3.15-}, we obtain
$\widehat\g^{\prime}=- \frac{1}{\sqrt{2}}\big(E\Phi_1^{\prime}\big)(x)^*\b(x)$. Hence, the third equality in \eqref{A8}, namely
$\widehat\g^{\prime}J\widehat\g^*=0$, follows.
\end{proof}
Recall that, according to Remark  \ref{RkG}, $\g(x)$ satisfies the relations
\begin{align}\label{A8'} &
\g(0)=\frac{1}{\sqrt{2}}\begin{bmatrix}-I_p & I_p \end{bmatrix}, \quad \g^{\prime}(x)J\g(x)^* \equiv 0,
\end{align}
from which $\g J \g^*\equiv -I_p$ follows. In view of \eqref{1.21} and \eqref{1.23+}, the matrix function $\b(x)$
satisfies the relations
\begin{align}\label{A20-} &
\b(0)=\frac{1}{\sqrt{2}}\begin{bmatrix}I_p & I_p \end{bmatrix}, \quad \b^{\prime}(x)J\b(x)^* \equiv 0, \quad \b(x)J\g(x)^*\equiv 0,
\end{align}
and $\b J\b^*\equiv I_p$ follows. The recovery of  $\b(x)$ is described in   a useful and somewhat more general proposition below (compare with \cite[Proposition~2.2]{ALS-HLB}).
\begin{Pn} \label{PnBG} 
Let a given $p\times 2p$ matrix function  $\g(x)$ $(0\leq x \leq 2\ell)$ have absolutely continuous entries
and satisfy relations \eqref{A8'}. Then, there is a unique matrix function  $\b(x)$ on $[0,2\ell]$ $($with absolutely continuous entries$)$
such that \eqref{A20-} holds.
This $\b$ is recovered from $\g$ using the equality 
\begin{align}\label{A20} &
\b(x)=\chi(x)\wh\b(x), \quad \wh \b(x):=\begin{bmatrix}(\g_2(x)^*)^{-1} & - (\g_1(x)^*)^{-1}\end{bmatrix}, 
\end{align}
where the  $p\times p$ matrix function $\chi(x)$ is the unique solution of the matrix differential equation
\begin{align}\label{A21} &
\chi^{\prime}(x)=-\chi(x)\wh\b^{\prime}(x)J\wh \b(x)^*(\wh \b(x) J \wh \b(x)^*)^{-1}
\end{align}
with the initial condition $\chi(0)=(1/2) I_p$.
\end{Pn}
\begin{proof} Since $\g J \g^*\equiv -I_p$, $\g_1$ and $\g_2$ are invertible and the matrix function $\wh \b(x)$
is well defined in \eqref{A20}. Moreover, we have $\wh \b(x)J\g(x)^*\equiv 0$.  
Thus, the class of  absolutely  continuous $\b$ satisfying $\b J \g^*\equiv 0$ (i.e., the third equality in \eqref{A20-}) is described by formula \eqref{A20}, where $\chi$  are absolutely continuous
 $p\times p$ matrix functions.   (Note that if $\b J \b^*\equiv I_p$ holds, and so $\b(x)$ is non-degenerate, $\chi(x)$ in \eqref{A20} should be non-degenerate as well.)
 
 Clearly, in our case of $\g(0)$ given in \eqref{A8'},
 $\chi(0)=(1/2) I_p$ is the  unique initial condition such that the first equality in \eqref{A20-} holds.

 Finally,   the equality $\g J\g^*\equiv -I_p$ implies the invertibility of  $\wh \b(x) J \wh \b(x)^*$. Hence, taking into account
 \eqref{A20}, we see that \eqref{A21} is equivalent to the identity $\b^{\prime}J \wh \b^*\equiv 0$.  
 
 Thus $\b$, which we recover using \eqref{A20}, \eqref{A21} and the initial condition
 $\chi(0)=(1/2) I_p$, satisfies \eqref{A20-}.
  The uniqueness of  $\b$ (or, equivalently, $\chi$) follows from the fact that  $\chi$ should be non-degenerate
 (see above) and so $\b^{\prime}J  \b^*\equiv 0$ yields $\b^{\prime}J \wh \b^*\equiv 0$, which implies \eqref{A21}.
  \end{proof}

{\bf Acknowledgments}  {This research    was supported by the
Austrian Science Fund (FWF) under Grant  No. Y-963.} 

\end{document}